\documentclass[11pt,reqno]{amsart}
\usepackage{fullpage} % Package to use full page
\usepackage{parskip} % Package to tweak paragraph skipping
\usepackage{amsmath,amssymb,amsfonts,amsthm}
\usepackage{mathrsfs}
\usepackage{graphicx,textcomp}

\usepackage[cmtip,all]{xy}

\usepackage{TqcCommands}%ZW package
\usepackage{tikz, braids} % For graphics, drawing braids

\newtheorem{theorem}{Theorem}[section]

\theoremstyle{definition}
\newtheorem{definition}{Definition}[section]
\newtheorem{example}{Example}[section]

\theoremstyle{remark}
\newtheorem{remark}{Remark}[section]

\theoremstyle{proposition}
\newtheorem{prop}{Proposition}[section]

\theoremstyle{conjecture}
\newtheorem{conjecture}{Conjecture}[section]

\numberwithin{equation}{section}
\DeclareMathOperator{\Aut}{Aut}
\DeclareMathOperator{\End}{End}
\DeclareMathOperator{\FSexp}{FSexp}
\DeclareMathOperator{\sVec}{sVec}
\DeclareMathOperator{\tr}{tr}
\DeclareMathOperator{\vect}{Vec}
\newcommand{\uPic}{\underline{\mathrm{Pic}}}

\begin{document}

\title{Mathematics of Topological Quantum Computing}


%    author one information
\author{Eric C. Rowell}
\address{Department of Mathematics\\
    Texas A\&M University \\
    College Station, TX 77843\\
    U.S.A.}
\curraddr{}
\email{rowell@math.tamu.edu}
\thanks{The authors thank Susan Friedlander for soliciting this survey, the anonymous referees for many helpful comments and suggestions, and Xingshan Cui and JM Landsberg for the same.  Rowell is partially supported by NSF grant DMS-1410144, and Wang by NSF grant DMS 1411212.}

%    author two information
\author{Zhenghan Wang}
\address{Microsoft Station Q and Dept of Mathematics\\
    University of California\\
    Santa Barbara, CA 93106-6105\\
    U.S.A.}
\curraddr{}
\email{zhenghwa@microsoft.com\\zhenghwa@math.ucsb.edu}
\thanks{}

%    \subjclass is required.

\date{}

\dedicatory{}

%    Abstract is required.
\begin{abstract}

In topological quantum computing, information is encoded in \lq\lq knotted" quantum states of topological phases of matter, thus being locked into topology to prevent decay.  Topological precision has been confirmed in quantum Hall liquids by experiments to an accuracy of $10^{-10}$, and harnessed to stabilize quantum memory.   In this survey,  we discuss the conceptual development of this interdisciplinary field at the juncture of mathematics, physics and computer science.  Our focus is on computing and physical motivations, basic mathematical notions and results, open problems and future directions related to and/or inspired by topological quantum computing.

\end{abstract}

\maketitle


\section{Introduction}

When the \lq\lq decision problem" was initially posed by Hilbert in 1928, he presupposed a definition of an algorithm.   But a universally accepted formalization occurred only in Turing's 1936 landmark paper \lq\lq On Computable Numbers, with an Application to the Entscheidungsproblem" \cite{Tu38}.  As described in the biography of Turing \cite{Hodges}, one afternoon during the early summer of 1935, Turing, while lying in a meadow at Grantchester in London after a long run, had the inspiration to abstract a human being calculating with pencil and paper into \lq\lq a mechanical process"---now bearing his name: a Turing machine.  Subsequently, all algorithms from Euclidean to AlphaGo share the same mathematical abstraction.  Even Shor's quantum factoring algorithm can be run on a Turing machine, though a slow one.  Turing understood quantum physics well.  It begs the question why Turing did not pursue a quantum generalization of his machines given the intellectual tools he had at hand. 

Turing unified data and program by his universal machine in \cite{Tu38}.  Since then, computability became part of mathematics following independent work of Turing, Church, and Post.   That a real number is computable should be comparable to a number being algebraic.  The absoluteness of computability is enshrined into the Church or Church-Turing thesis: anything computable by a discrete-state-machine with finite means is computable by a Turing machine.  No serious challenge to the Church thesis has ever appeared\footnote{Hypercomputation is not considered as a viable possibility.}.  But in 1994, Shor's quantum algorithm for factoring integers \cite{Sh94} posed a serious challenge, not to the Church thesis, but to the polynomial extension of the Church thesis: anything efficiently computable can be efficiently computed by a Turing machine.  While there seems to be a unique notion of computability, it is likely that the notion of {\it efficient} computability would diverge.

Definiteness is an important feature of the classical world.  Thus, computing tasks can be formalized as Boolean functions.  Definability is not the same as computability because there exist plenty of non-computable numbers.  As currently formulated, quantum computing will not impact computability.  But classical notions can be smeared into quantum weirdness such as Schr\"odinger's cat.  The promise of quantum computing is a vast leap in the speed of processing classical information using quantum media.

The unit of information is bit.  Qubits are the quantum embodiment of bits using the two characteristic attributes of quantum mechanics---superposition and entanglement.  Qubits are computationally powerful, but notoriously fragile---the outside world is constantly \lq\lq looking at" (measuring) the quantum system, which leads to the decoherence of the quantum states.

Mathematically, a qubit is the abstraction of all quantum systems with two-dimensional Hilbert spaces of states.  We call the two-dimensional Hilbert space $\mbbC^2$ with preferred basis $|0\rangle$ and $|1\rangle$ {\it a qubit}.  Therefore, a qubit, utilizing superposition, can have any non-zero vector $|\psi\rangle \in \mbbC^2$ as a state, but computability considerations will constrain us only to states that can be algorithmically prepared.  Experiment will further force us to work with only finite precision states.  Einstein's \lq\lq spooky action at a distance"---quantum entanglement---is realized in states of multi-qubits ${(\mbbC^2)}^{\otimes n}$, where a state $|\psi\rangle$ is entangled if $|\psi\rangle$ cannot be written as a tensor product.

Anyons, generalizations of bosons and fermions in flatland, are topological quantum fields materialized as finite energy particle-like excitations in topological phases of matter, the subject of the 2016 Nobel and Buckley prizes in physics.  Like particles, they can be moved, but cannot be created or destroyed locally.  There are two equivalent ways to model anyon systems.  We can focus on the ground state manifold $V(Y)$ of an anyonic system on any possible space $Y$, and then the anyon system is modeled in low energy by a unitary $(2+1)$-TQFT.  An alternative is to consider the fusion and braiding structures of all elementary excitations in the plane.  The anyon system is then equivalently modeled by a unitary modular category.  The two notions unitary $(2+1)$-TQFT and unitary modular category are essentially the same \cite{Tu94}.  Therefore, anyon systems can be modeled either by unitary TQFTs or unitary modular categories.  In the modular category model, an anyon $X$ is a simple object that abstracts an irreducible representation of some algebra of symmetries.

Topological quantum computing (TQC) solves the fragility of the qubits at the hardware level\footnote{The software solution \cite{Shor95}, fault-tolerant quantum computation, requires daunting overhead for practical implementation.} by using topological invariants of quantum systems.  Information is encoded non-locally into topological invariants that spread into local quantities just as the Euler characteristic spreads into local curvature by the Gauss-Bonnet theorem.  Nature does provide such topological invariants in topological phases of matter such as the FQH liquids and topological insulators.  The topological invariant for TQC is the ground state degeneracy in topologically ordered states with non-abelian anyons.

Can the Jones polynomial $J(L,q)$ of oriented links $L$ evaluated at $q=e^{\frac{2\pi i}{r}}, r=1,2,3,\ldots$ be calculated by a computing machine quickly? The fact that the Jones evaluation is $\#P$-hard if $r\neq 1,2,3,4,6$ \cite{Ve91,JVW90,We93} is one of the origins of topological quantum computing \cite{Fr98}.  Freedman wrote in \cite{Fr98}:
\lq\lq Non-abelian topological quantum field theories exhibit the mathematical features necessary to support a model capable of solving all $\#P$ problems, a computationally intractable class, in polynomial time. Specifically, Witten \cite{Wi89} has identified expectation values in a certain SU(2) field theory with values of the Jones polynomial \cite{Jo85} that are $\#P$-hard \cite{JVW90}. This suggests that some physical system whose effective Lagrangian contains a non-abelian topological term might be manipulated to serve as an analog computer capable of solving NP or even $\#P$ hard problems in polynomial time. Defining such a system and addressing the accuracy issues inherent in preparation and measurement is a major unsolved problem."

Another inspiration for topological quantum computing is fault-tolerant quantum computation by anyons \cite{Ki03}.  Kitaev wrote in \cite{Ki06}:
\lq\lq A more practical reason to look for anyons is their potential use in quantum computing. In \cite{Ki03}, I suggested that topologically ordered states can serve as a physical analogue of error-correcting quantum codes.  Thus, anyonic systems can provide a realization of quantum memory that is {\it protected from decoherence}.  Some quantum gates can be implemented by braiding; this implementation is {\it exact} and does not require explicit error-correction.  Freedman, Larsen, and Wang \cite{FLW02} proved that for certain types of non-abelian anyons braiding enables one to perform universal quantum computation.  This scheme is usually referred to as {\it topological quantum computation}."

The theoretical resolution of Freedman's problem \cite{FKW02}, Kitaev's idea of inherently fault-tolerant quantum computation \cite{Ki03}, and the existence of universal anyons by braidings alone \cite{FLW02} ushered in topological quantum computing \cite{FKLW03}.

In TQC, information is encoded in \lq\lq knotted" quantum states, thus being locked into topology to prevent decay.  Topological precision has been confirmed in quantum Hall liquids by experiments to an accuracy of $10^{-10}$, and has been harnessed to stabilize quantum memory.  TQC has been driving an interaction between mathematics, physics, and computing science for the last two decades. In this survey, we will convey some of the excitement of this interdisciplinary field.

We need physical systems that harbor non-abelian anyons to build a topological quantum computer.  In 1991, Moore and Read, and Wen proposed that non-abelian anyons exist in certain fractional quantum Hall (FQH) liquids \cite{MR91,We91}.  Moore and Read used conformal blocks of conformal field theories (CFTs) to model fractional quantum Hall states, and Wen defined topologically ordered states, whose effective theories are topological quantum field theories (TQFTs).  Thus, $(2+1)$-TQFT and $(1+1)$-CFT form the foundations of TQC.  The algebraic input of a TQFT and topological properties of a CFT are encoded by a beautiful algebraic structure---a {\it modular tensor category}, which will be simply called a modular category.  The notion of a modular tensor category was invented by Moore and Seiberg using tensors \cite{MS89}, and its equivalent coordinate-free version modular category by Turaev \cite{Turaev92}.  Thus, modular categories, algebraically underpinning TQC, are algebraic models of anyon systems.

In 1999, Read and Rezayi \cite{RR99} suggested a connection between FQH liquids at filling fractions $\nu=2+\frac{k}{k+2}$ and $SU(2)_k$-Witten-Chern-Simons theories--mathematically Reshetikhin-Turaev TQFTs--for $k=1,2,3,4$.  In 2006, interferometer experiments of FQH liquids were proposed to show that the evaluations at $q=i$ of the Jones polynomial for certain links directly appear in the measurement of electrical current \cite{BKS06,SH06}.  In 2009, experimental data consistent with the prediction were published \cite{WPW09}.  The current most promising platform for TQC is nanowires and topological protection has been experimental confirmed \cite{Mou12,Marcus16}.  The next milestone will be the experimental confirmation of non-abelian fusion rules and braidings of Majorana zero modes \cite{lutchyn17,Scalable17}.

Progress towards building a useful quantum computer has accelerated in the last few years.  While it is hard to characterize our current computing power in terms of a number of qubits, it is clear that a working quantum computer with one hundred qubits would perform tasks that no classical computer can complete now.  Since TQC does not have a serious scaling issue, when one topological qubit is constructed, a powerful quantum computer is on the horizon.

There are many interesting open questions in TQC including the classification of mathematical models of topological phases of matter such as TQFTs, CFTs and modular categories, and the analysis of computational power of anyonic quantum computing models.  Classification of modular categories is achievable and interesting both in mathematics and physics \cite{RSW09,BNRW151,BNRW152}.  A recent result in this  direction is a proof of the rank-finiteness conjecture \cite{BNRW151}.  An interesting open question in the second direction is the property $F$ conjecture that the braid group representations afforded by a simple object $X$ in a modular category have finite images if and only if its squared quantum dimension $d_X^2$ is an integer \cite{NR11,ERW}.  Many open problems, for example the theories for fermions and three spatial dimensions, will be discussed in the survey.

The arenas of mathematics, computer science, and physics are Mind, Machine, and Nature.  Machine learning has come a long way since Turing's paper \cite{Turing50}.  AlphaGo is an example of the amazing power of Machines.  The question \lq\lq Can Machine think" is as fresh as now as it was in the 1950s.  We are at an important juncture to see how the three worlds would interact with each other.  TQC is the tip of an iceberg that blurs the three.  The authors' bet is on Nature, but we suspect that Nature has her eye on Machine.  An inevitable question that will soon confronts us is {\it how are we going to adapt when quantum computers become reality}?

In this survey, we focus on the mathematics of TQC as in \cite{Wa10}.  For the physical side, we recommend \cite{NSSFD08,Pr99, Pa12,Majorana15}.  Three fundamental notions for TQC are: modular category, $(2+1)$-TQFT, and topological phase of matter.  Though closely related, the definition of modular category is universally accepted while that of TQFT varies significantly.  For our applications, we emphasize two important principles from physics: locality and unitarity.  Locality follows from special relativity that nothing, including information, can travel faster than light, whereas unitarity is a requirement in quantum mechanics.  We define unitary $(2+1)$-TQFTs adapting the definitions of Walker and Turaev \cite{Walker91, Tu94}.  We propose a mathematical definition of 2D topological phases of matter using the Hamiltonian formalism in section 4.  The content of the survey is as follows:  In section 2, we give an introduction to TQC and mathematical models of anyons.  In section 3, we cover abstract quantum mechanics and quantum computing.  In section 4,  we first lay foundations for a mathematical study of topological phases of matter, then define $2D$ topological phases of matter.  In section 5, we analyze the computational power of anyonic computing models.  Section 6 is a survey on the structure and classification of modular categories.  In section 7, we discuss various extensions and open problems.  We conclude with two eccentric research directions in section 8.

The wide-ranging and expository nature of this survey makes it impossible to obtain any reasonable completeness for references.  So we mainly cite original references, expository surveys, and books.

\section{The ABC of Topological Quantum Computing}

We introduce three basic notions in topological quantum computing (TQC): {\underline{A}}nyons, {\underline{B}}raids, and {\underline{C}}ategories.  One salient feature of TQC is the extensive use of graphical calculus.  Space-time trajectories of anyons will be represented by braids, and more general quantum processes such as creation/annihilation and fusion by tangles and trivalent graphs.  Algebraically, anyon trajectories will be modeled by morphisms in certain unitary modular categories (UMCs). A UMC can be regarded as a computing system: the morphisms are circuits for computation.  After motivating the axioms of a modular category using anyon theory, we define this important notion and explain anyonic quantum computing.  Other elementary introductions to TQC include \cite{Wang06,Wang13,DRW16,Rowell16}.

\subsection{Anyons and Braids}

It is truly remarkable that all electrons, no matter where, when, and how they are found, are identical.  Elementary particles\footnote{Elementary particles are elementary excitations of the vacuum, which explains why they are identical.  In TQC, anyons are elementary excitations of some two dimensional quantum medium, so they are called quasi-particles.  The vacuum is also a very complicated quantum medium, so we will simply refer to anyons as particles.  This use of particle sometimes causes confusion since then there are two kinds of particles in an anyonic physical system: anyons, and the constituent particles of the quantum medium such as electrons in the fractional quantum Hall liquids where the anyons emerge.} are divided into bosons and fermions.  Consider $n$ quantum particles $X_i, 1\leq i \leq n$, in $\mbbR^3$ at distinct locations $r_i$, then their quantum state is given by a wave function\footnote{We ignore other degrees of freedom such as spins.} $\Psi(r_1,...,r_i,...,r_j,...,r_n)$.  If we exchange $X_i$ and $X_j$ along some path so that no two particles collide during the exchange, then $\Psi(r_1,...,r_j,...,r_i,...,r_n)=\theta\cdot \Psi(r_1,...,r_i,...,r_j,...,r_n)$ for some complex number $\theta$.  If we repeat the same exchange, the $n$ particles will return back exactly, hence $\theta^2=1$.  Particles with $\theta=1$ are bosons, and $\theta=-1$ fermions.  It follows that if we perform any permutation $\sigma$ of the $n$ particles, then $\Psi(r_{\sigma(1)},...,r_{\sigma(n)})=\pi(\sigma)\cdot \Psi(r_1,...,r_n)$, where $\pi(\sigma)$ is the sign of the permutation $\sigma$ if the particles are fermions.  Therefore, the statistics of elementary particles in $\mbbR^3$ is a representation of the permutation group $\fS_n$ to $\mbbZ_2=\{\pm 1\}$.

Abstractly, the positions of $n$ identical particles living in a space $X$ form the configuration space $C_n(X)$ and the quantum states of such $n$ particles form some Hilbert space $V_n$.  When $n$ particles move from one configuration in a state $v_0\in V_n$ back to the initial configuration, the trajectory is a loop in $C_n(X)$ and induces a change of states from $v_0$ to some potentially different state $v_1\in V_n$.  If the particles are topological, then the change of states depends only on the homotopy class of loops, hence gives rise to an action of the fundamental group of $C_n(X)$ on $V_n$.  The collection of representations $\{V_n\}$ for all $n$ is called the statistics of the particle.

Now we examine two implicit facts that are used in the above discussion.
Firstly, we assume that the final state of the $n$ particles does not depend on the exchange paths. This is because any two paths for the exchange are isotopic\footnote{This is the fact there are no knotted simple loops in $\mbbR^4$.}.  But technology has made it realistic to consider particles confined in a plane.  Then isotopic classes of paths for exchange form the braid groups $\mcB_n$ instead of the permutation groups $\fS_n$.  Repeating the argument above, we can only conclude that $\theta$ is on the unit circle $\mathbb{U}(1)$.  Potentially $\theta$ could be any phase $\theta\in \mathbb{U}(1)$.  Particles with {\underline{any}} phase $\theta$ are dubbed {\underline{any}}ons by Wilczek\footnote{Stability consideration restricts possible phases to $\theta=e^{2\pi i s}$, where $s$ is a rational number. The value $s$ is related to the strength of a fictitious flux, so an irrational value would mean infinite precision of a magnetic field, which is not stable.}.  For the history and classical references on anyons, see \cite{wil90}. Now they are referred to as abelian anyons because their statistics are representations of the braid groups $\mcB_n$ into the abelian group $\mathbb{U}(1)$.

The second fact used is that there is a unique state $\Psi(r_1,...,r_n)$ when the positions of the $n$ particles are fixed.  What happens if there is more than one linearly independent state? If the $n$ particles start with some state, they might come back to a superposition after a time evolution.  Suppose $\{e_i\}, 1\leq i \leq m$, is an orthonormal basis of all possible states, then starting in the state $e_i$, the $n$ particles will return to a state $\sum^m_{j=1}U_{ij}e_j$, where $(U_{ij})_{1\leq i,j\leq m}$ is an $m\times m$ unitary matrix.  The statistics of such particles could be high dimensional representations of the braid group\footnote{We might wonder why there cannot be elementary particles that realize higher dimensional representations of the permutation groups. Such statistics is called parastatistics.  It has been argued that parastatistics does not lead to fundamentally new physics as non-abelian statistics do.} $\mcB_n\rightarrow \mathbb{U}(m)$.  When $m>1$, such anyons are non-abelian.

Space-time trajectories of $n$ anyons form the \emph{n-strand braid group} $\mcB_n$.  Mathematically, $\mcB_n$ is the the motion group of $n$ points in the disk $D^2$ given by the presentation

\begin{align*}\mcB_n =  \langle \sigma_1, \sigma_2, \ldots, \sigma_{n-1} \mid &\sigma_i\sigma_j=\sigma_j \sigma_i \text{ for } |i -j| \ge 2,\\ &\sigma_i\sigma_{i+1}\sigma_i=\sigma_{i+1}\sigma_i\sigma_{i+1}, \text{ for }  1\leq i,j \leq n-1  \rangle\end{align*}

The first relation is referred to as far commutativity and the second is the famous braid relation. The terminologies become clear when one considers the graphical representation of elements of the braid group.  A fun example is the following $4$-strand braid $b=\sigma_3^{-1}\sigma_2^2 \sigma_{3}^{-1}\sigma_1^{-1}$, once drawn by Gauss.

$$\begin{tikzpicture}[scale=.5]
 \braid[number of strands=4] (braid)  a_1 a_3 a_2^{-1} a_2^{-1}  a_3 ;

\end{tikzpicture} $$

We have an exact sequence of groups $$1\rightarrow P_n \rightarrow \mcB_n\rightarrow \fS_n\rightarrow 1,$$ where $P_n$ is called the $n$-strand pure braid group.  Note that from this sequence, any representation of $\fS_n$ leads to a representation of $\mcB_n$.

\subsection{Tensor Categories and Quantum Physics}

Tensor categories, categorifications of rings, are linear monoidal categories with two bifunctors $\oplus, \otimes$, which are the categorifications of sum $+$ and multiplication $\times$ of the ring.  The strict associativity of the multiplication $\times$ is relaxed to some natural isomorphisms of tensor product $\otimes$.  A standard reference on tensor categories is \cite{etingof15}.

There is a philosophical explanation for why tensor category theory is suitable for describing some quantum physics.\footnote{During the conference Topology in Condensed Matter Physics at American Institute of Mathematics in 2003, Xiao-Gang Wen asked the second author what should be the right mathematical framework to describe his topological order.  \lq\lq Tensor Category Theory" was his reply.}  In quantum physics, we face the challenge that we cannot \lq\lq see" what is happening.  So we appeal to measurements and build our understanding from the responses to measuring devices.  Therefore, quantum particles (e.g. anyons) are only defined by how they interact with other particles, and their responses to measuring devices.  In tensor category theory, objects are usually the unknown that we are interested in.  An object $X$ is determined by the vector spaces of morphisms $\Hom(X,Y)$ for all $Y$ in the tensor category.  Therefore, it is natural to treat objects in a tensor category as certain special quantum states such as anyons, and the morphisms as models of the quantum processes between them.\footnote{The tensor product is needed to model many anyons at different space locations.}  The UMC model of anyon systems is such an example \cite{FKLW03,Ki06,Wa10}.

Another philosophical relation between tensor categories and quantum physics comes from a similarity between quantization and categorification.  A famous quote by Nelson is: \lq\lq first quantization is a mystery, but second quantization is a functor".  The quantization of a finite set $S$ is a good illustration.  First quantization is the process going from a classical system to a quantum system that is modeled by a Hilbert space of quantum states with a Hamiltonian.  In the case of a finite set $S$ as the classical configuration space of a single particle, then quantization is simply the linearization of $S$---the Hilbert space is just $\mbbC[S]$ with basis $S$.  Second quantization is the process going from a single particle Hilbert space to a multi-particle Hilbert space.  For simplicity, consider a fermion with a single particle Hilbert space $V$ of dimension$=n$, then the multi-particle Fock space is just the exterior algebra $\wedge^{*}V$ of dimension $2^n$.  Hence, second quantization is the functor from $V$ to $\wedge^{*}V$.
The process of de-quantization is measurement:  when we measure a physical observable $\mcO$ at state $|\Psi\rangle$, we arrive at a normalized eigenvector $e_i$ of $\mcO$ with probability $p_i=\langle e_i|\mcO|\Psi\rangle$.  A basis consisting of eigenstates of an observable $\mcO$ are the possible states that we obtain after measuring $\mcO$.  Therefore, eigenstates of observables are the {\it reality} that we \lq\lq see" classically.

According to Kapranov and Voevodsky, the main principle in category theory is: \lq\lq In any category it is unnatural
and undesirable to speak about equality of two objects".  The general idea of categorification is to weaken the nonphysical notion of equality to some natural isomorphism.  Naively, we want to replace a natural number $n$ with a vector space of dimension $n$.  Then the categorification of a finite set $S$ of $n$ elements would be $\mbbC[S]$ of dimension $n$.  It follows that the equality of two sets $S_1,S_2$ should be relaxed to an isomorphism of the two vector spaces $\mbbC[S_1]$ and $\mbbC[S_2]$.  Isomorphisms between vector spaces may or may not be functorial as the isomorphisms between $V$ and $V^{**}$ or $V^{*}$ demonstrate.  It follows that categories are more physical notions than sets because the ability to instantaneously distinguish two elements of any set is nonphysical\footnote{We regard both infinite precision and infinite energy as nonphysical.  If there is a minimal fixed amount of energy cost to distinguish any two elements in a set, it would cost an infinite amount of energy to distinguish instantaneously every pair of elements of an infinite set.}.

\subsection{Gedanken Experiments of Anyons}

A topological phase of matter (TPM) is an equivalence class of lattice Hamiltonians or ground states which realizes a unitary topological quantum field theory (TQFT) at low energy (a more detailed discussion on lattice Hamiltonians and TPMs is in chapter 4).  Anyons are elementary excitations in TPMs, therefore they are special point-like quantum states which form a closed particle system.
How can we model a universe of anyons following the laws of quantum physics, locality, and topological invariance?  

An anyon system is much like the world of photons and electrons in quantum electrodynamics.  First there are only finitely many anyon types in each anyon universe, which form the label set $L=\{0,\ldots,r-1\}$ of \lq\lq labels"\footnote{Other names include super-selection sectors, \lq\lq anyon charges", \lq\lq anyon types", and \lq\lq topological charges".}.  Let $A_L=\{X_i\}_{i\in L}$ be a representative set of anyon types: one for each anyon type.  The ground state or \lq\lq vacuum" $X_0=\unit$ of the topological phase is always included as an invisible anyon labeled by $0$.  It is important to distinguish anyons and their types, so we will usually use upper case letters for physical anyons and the corresponding lower cases for their types.  An anyon $X$ with type $x$ has an anti-particle $X^*$ with type denoted by $\hat{x}$. The ground state $\unit$ is its own antiparticle so $\unit^*\cong\unit$ and $\hat{0}=0$.  The number of anyon types $r$ is called the rank of the anyon system.

\subsubsection{Ground state degeneracy and quantum dimension}
Consider an oriented closed two dimensional space $Y$ with some anyons $\{\tau_i\}$ residing at $\{q_i\}$.  Topological invariance means that the local degrees of freedom such as positions $\{q_i\}$ will not affect universal physical properties.  Even though anyons are point-like, they do occupy some physical space.  There are two equivalent ways to visualize anyons: either imagine an anyon as a colored puncture in $Y$ with a signed infinitesimal tangent vector, or as an infinitesimal disk whose boundary circle has a sign, a based point and a color.  By a color, we mean an anyon, not its type.  We will explain the sign and base point later in section 2.6.1.  When anyons $\{\tau_i\}$ are fixed at positions $\{q_i\}$, all quantum states of the topological phase form a Hilbert space $\mcL(Y;\tau_i,q_i,p_i)$, which also depends on other local degrees of freedom $p_i$ such as momenta and spins.  For topological phases, we are only interested in the topological degrees of freedom $V(Y;\tau_i)$ in the low energy states because high energy states easily decay (more explanation is in chapter 4).  The well-defined low energy topological degrees of freedom due to the energy gap form the relative ground state manifold\footnote{Physicists use manifold here to mean multi-fold, but the ground state Hilbert space $V(Y;\tau_i)$ here or later $\mcL_{\lambda_0}$ is also the simple manifold $\mathbb{R}^n$ for some $n$.} $V(Y;\tau_i)$.

Locality implies that a state in $V(Y;\tau_i)$ can be constructed from states on local patches of $Y$.  Every compact surface $Y$ has a DAP decomposition:  a decomposition of $Y$ into disks, annuli, and pairs of pants.  The new cutting circles of a DAP decomposition should have some boundary conditions $\ell$ that allow the reconstruction of the original quantum state from ground states on disks, annuli, and pairs of pants.  A physical assumption is that the stable boundary conditions $\ell$ of a cutting circle are in one-one correspondence with anyons in the complete representative set $A_L=\{X_i\}_{i\in L}$, which is a version of bulk-edge correspondence.  Then locality is encoded in a gluing formula of the form:

$$V(Y;\tau_i)\cong \bigoplus_{X_i\in A_L}V(Y_{\textrm{cuts}}; \tau_i, X_i,{X_i}^*).$$

Therefore, general anyon states on any space $Y$ can be reconstructed from anyon states on disks, annuli, and pairs of pants together with some general principles. The topology of the disk, annulus, and pair of pants detect the vacuum, anti-particles, and fusion rules. 

An important quantum number of an anyon $X$ is its quantum dimension $d_X$---a positive real number $\geq 1$.  The quantum dimension $d_X$ of an anyon $X$ determines the asymptotic growth rate when $n$ identical anyons $X$ are confined to the sphere: the dimension of the degeneracy ground state manifold $V_{X,n,\unit}=\textrm{Hom}(X^{\otimes n}, \unit)$ grows as $(d_X)^n$ as $n\rightarrow \infty$.  Therefore, an anyon $X$ leads to degeneracy in the plane if and only if its quantum dimension $d_X >1$.
An anyon $X$ is {non-abelian} if $d_X>1$, otherwise $d_X=1$ and is abelian.

Quantum dimensions of anyons are related to topological entanglement entropy, which is an important way to characterize topological phases.  Mathematically, they are the same as the Frobenius-Perron dimensions in section 2.4.3.  Quantum dimensions have interesting number theoretical properties, for example they are $S$-units in a cyclotomic field \cite{BNRW151}.

\subsubsection{Vacuum and anti-particles}

Since a disk $D^2$ has trivial topology, i.e. is contractible, it cannot support any nontrivial topological state other than the unique topological ground state.  This disk axiom $\dim(V(D^2;X_i))=\delta_{0i}$ essentially implies that topological ground state manifolds $V(Y;\tau_i)$ are error-correcting codes\footnote{When a TQFT is realized by a lattice model with a microscopic Hamiltonian, local operators supported on sufficiently small disks cannot act on the ground states other than as a scalar multiple of the identity due to the disk axiom.  This property is a characterization of a quantum error correcting code.}: $V(Y;\tau_i)\subset \mcL(Y;\tau_i,q_i,p_i)$.  Since an anyon cannot change its topological charge across the annular region between the two boundary circles, a pair of anyons coloring the boundary circles cannot support a nonzero topological state on the annulus $\mcA$ unless they are dual particles to each other (dual comes from the opposite induced orientation): $\textrm{dim}(V(\mathcal{A};X_i,X_j))=\delta_{i\hat{j}}$.

\subsubsection{Creation/annihilation, fusion/splitting, and fusion rules}

\begin{figure}[h]\includegraphics[width=2.5in]{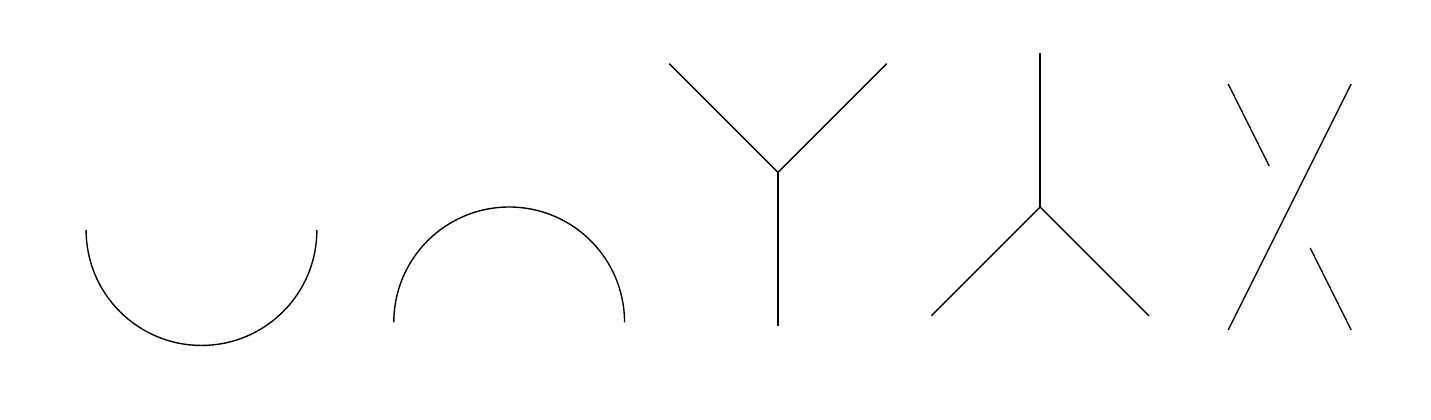}
\caption{\label{fig6}Elementary Events.}
\end{figure}

Fix an anyon system with anyon representatives $A_L=\{X_i\}_{i\in L}$.  Suppose there are $n$ anyons $X_i, i=1,2,...,n$ in the interior of the disk $D^2$ and the outermost boundary of the disk is colored by some anyon $X_\infty\in A_L$.  The dimension of $V(D^2;X_i, X_\infty)$ can be inductively found by the gluing formula.  If there are no anyons in the disk, by the disk axiom, then there will be no quantum states in the disk unless $X_\infty=\unit$.  Choosing a normalized state $|0\rangle$ in $V(D^2;\unit)$ as the no-anyon state---the vacuum for the disk, then we can create pairs of anyons from  $|0\rangle$.  It is a fundamental property of anyons that a single anyon cannot be created from the ground state---a conservation law.  How would we measure anyon states?  The most important one is by fusion: we bring two anyons $X_i$ and $X_j$ together and see what other anyons $X_k$ would result from their fusion.  This process is represented by a {\large $\mathsf{Y}$} if time goes from top to bottom.  If time goes up, then the same picture represents the splitting of a single anyon into two.  Using our picture of an anyon as a small disk, we know {\large $\mathsf{Y}$} really should be thickened to a pairs of pants.  If we label the three circles with all possible anyons $X_i\in A_L$, we obtain a collection of non-negative integers $N_{ij}^k=\dim{V(D^2;X_i,X_j,X_k)}$.  The collection of integers $N_{ij}^k$ is called the fusion rule. Fusion rules are also written as $x_i\otimes x_j=\bigoplus_{k}N_{ij}^k x_k.$
Therefore, the elementary events for anyons are creation/annihilation, fusion/splitting, and braiding of anyons, see Figure \ref{fig6}.

The quantum dimension of an anyon can be easily computed from its fusion rules: regard all anyon types as unknown variables and the fusion rules as polynomial equations, then the maximal real solutions of these polynomial equations are the quantum dimensions.

\subsubsection{Non-abelian statistics}

Non-abelian statistics is a fundamentally new form of particle interactions.  This \lq\lq spooky action" is a manifestation of the entanglement of the degenerate ground states in  $V(Y;\tau_i)$.  Besides its general interest as a new form of particle interaction, non-abelian statistics underlies the idea of TQC---the braiding matrices are inherently fault-tolerant quantum circuits \cite{Ki03,Fr98, FKLW03}.
The ground state manifold $V(D^2;X,...,X,X_\infty)$ is a representation of the $n$-stand braid group $\mcB_n$. 
The {\it statistics of an anyon} $X$ are these representations: $\rho_{X,n}: \mcB_n \rightarrow \mathbb{U}(V(D^2;X,...,X,X_\infty))$.  The anyon $X$ is non-abelian if not all braid images form abelian subgroups.

\subsection{Mathematical Models of Anyons}

The first mathematical model of anyon systems is through unitary topological modular functors (UTMF) underlying the Jones polynomial \cite{FKLW03}.  We adapt the definition of modular functor as axiomatized by Turaev\footnote{Our modular functor is what Turaev called a weak rational 2D modular functor.} \cite{Tu94}, which is essentially the topological version of a modular category.  Therefore, an abstract anyon system can be either identified as a UTMF topologically or a UMC algebraically.  There is a subtle Frobenius-Schur indicator that complicates the axiomatization of a modular functor.  In the context of a modular category, it is then important to distinguish between an anyon and its equivalence class.  The distinction between an anyon and its type is emphasized in \cite{Wa10}.  Two anyons are equivalent in physics if they differ by local operators.  In a UMC, local operators correspond to the morphisms that identify the two simple objects. That an anyon is a simple object in a UMC seems to be first stated explicitly in \cite{Wa10}.

\subsubsection{Unitary topological modular functor}\label{Sec:UTMF}

The universe of an anyon system can be formalized into a unitary strict fusion category first.  Anyons and their composites will be objects of a category.  A composite of several anyons should be considered as their product, called a tensor product $\otimes$.  The ground state is the tensor unit $\unit$ as we can attach it to any anyon without changing the topological properties.  Therefore, we have a special monoidal class in the sense of Turaev\footnote{Our definition differs from Turaev's in two important aspects: first our input data is a strict fusion category rather than a monoidal class.  This stronger assumption is one way to implement Turaev's duality axiom 1.5.8 on page 245 in \cite{Tu94}. Secondly we add an algebraic axiom to ensure all mapping class group representation matrices can be defined within some number field, which is necessary for applications to quantum computing.} \cite{Tu94}.  Another important operation is a sum $\oplus$, which results from the fusion of anyons.

Recall from \cite{ENO1}, a \textit{fusion category}  $\mcC$ over $\mbbC$ is an abelian $\mbbC$-linear semisimple rigid monoidal category with a simple unit object $\unit$,
finite-dimensional morphism spaces and finitely many isomorphism classes of simple objects.  A fusion category is unitary if all morphism spaces are Hilbert spaces with certain properties \cite{Tu94,Wa10}.

Let $\Pi_\mcC$ be the set of isomorphism classes of simple objects of the fusion category $\mcC$. The set $\Pi_\mcC$ is called a label set in TQFTs, and the set of anyon types or topological charges in anyonic models.  The {\textit{rank}} of $\mcC$ is the finite number $r=|\Pi_\mcC|$, and we denote the members of $\Pi_\mcC$ by $\{0,\ldots,r-1\}$. We simply write $V_i$ for an object in the isomorphism class $i \in \Pi_\mcC$. By convention, the isomorphism class of $\unit$ corresponds to $0\in\Pi_\mcC$. The rigidity of $\mcC$ defines an involution $i \leftrightarrow \hat{i}$ on $\Pi_\mcC$ which is given by $V_{\hat{i}} \cong V_i^{*}$ for all $i \in \Pi_\mcC$.

For our definition of a UTMF (Definition 2.1), we define a projective version $\Hil$ of the category of finite dimensional Hilbert spaces.  The objects of $\Hil$ are finite dimensional Hilbert spaces and morphisms are unitary maps up to phases, i.e. equivalence classes of morphisms that differ by only a phase. A projective functor $F: \mcC \rightarrow \mcD$ is a map that satisfies all properties of a functor except that $F$ preserves composition of morphisms only up to a phase, i.e. $F(fg)=\xi(f,g) F(f) F(g)$ for some $\xi(f,g) \in \mathbb{U}(1)$ without any other conditions on $\xi(f,g)$.

Given a strict fusion category $\mcC$, we define the category $\Bord$ of colored oriented compact surfaces with boundary conditions $\mcC$: the objects of $\Bord$ are compact oriented surfaces $Y$ with oriented, based, colored boundary circles.  The empty set $\emptyset$ is considered as a surface with a unique orientation.  The boundary $\partial Y$ of $Y$ consists of oriented circles\footnote{The independent orientation may or may not agree with the induced orientation from the surface, which is the same as the sign in anyon theory.} with a base point, an orientation independent of the induced orientations from $Y$, and a simple object of $\mcC$---a color.  Morphisms of $\Bord$ are data preserving diffeomorphisms up to data preserving isotopy: a diffeomorphism which preserves the orientation of $Y$, the based points, orientations, and colors of boundary circles.

In the following, $\simeq$ denotes isomorphisms which are not necessarily canonical, while $\cong$ are functorial isomorphisms.  For a set of labels $\ell=\{l_i\}$, $X_\ell$ denotes the tensor product $\otimes_{l_i}X_{l_i}$ of representatives for each label $l_i$, which is well-defined as our category is strict.  For a strict fusion category $\mcC$, $A_{\mcC}=\{X_i\}_{i\in \Pi_{\mcC}}$ is a complete set of representatives of simple objects of $\mcC$.

\begin{definition}

Given a strict fusion category $\mcC$, a unitary topological modular functor $V$ associated with the strict unitary fusion category $\mcC$ is a symmetric monoidal projective functor\footnote{The monoidal and braiding structures are not projective.} from $\Bord \rightarrow \Hil$ satisfying the following additional axioms:

\begin{enumerate}

\item \textrm{Empty surface axiom:}
$V(\emptyset)=\mbbC.$

\item \textrm{Disk axiom:}

$V(D^2;X_i)\cong \begin{cases}
\mbbC & i=0,\\
0 & \textrm{otherwise}, \end{cases}$
\quad \quad where $D^2$ is a $2$-disk.

\item \textrm{Annular axiom:}

\[V(\mcA;X_i,X_j)\simeq \begin{cases} \mbbC & \textrm{if $i=\hat{j}$},\\ 0 & \textrm{otherwise}, \end{cases}\]
where $\mcA$ is an annulus and $i,j\in \Pi_{\mcC}$.  Furthermore, $V(\mcA;X_i,X_j)\cong \mbbC$ if  $X_i\cong X_j^{*}$.

\item \textrm{Disjoint union axiom:}

$V(Y_1 \sqcup Y_2;X_{\ell_1}\sqcup X_{\ell_2})\cong V(Y_1;X_{\ell_1})\otimes
V(Y_2;X_{\ell_2})$.
\\
The isomorphisms are associative, and
compatible with the mapping class group projective actions $V(f): V(Y)\rightarrow V(Y)$ for $f: Y\rightarrow Y$.

\item \textrm{Duality axiom:}

$V(-Y;X_{\ell})\cong V(Y;X_{\ell})^*$, where $-Y$ is $Y$ with the opposite orientation.
The isomorphisms are compatible with mapping class group projective actions, orientation reversal, and the disjoint union axiom as follows:

(i): The isomorphisms $V(Y) \rightarrow V(-Y)^*$ and $V(-Y)\rightarrow V(Y)^*$ are mutually adjoint.

(ii): Given $f: (Y_1;X_{\ell_1})\rightarrow (Y_2;X_{\ell_2})$ let $\bar{f}:
(-Y_1; X_{\ell_1}^*)\rightarrow (-Y_2;X_{\ell_2}^*)$ be the induced reversed orientation map, we have
$\langle x,y\rangle =\langle V(f)x,V(\bar{f})y\rangle $, where $x\in V(Y_1;X_{\ell_1})$, $y\in V(-Y_1;X_{\ell_1}^*)$.

(iii): $\langle \alpha_1\otimes \alpha_2, \beta_1\otimes \beta_2\rangle =\langle \alpha_1, \beta_1\rangle \langle \alpha_2,\beta_2\rangle$ whenever
\begin{align*}
\alpha_1\otimes \alpha_2 &\in V(Y_1\sqcup Y_2)\cong V(Y_1)\otimes V(Y_2),\\
\beta_1\otimes \beta_2 &\in V(-Y_1\sqcup -Y_2)\cong V(-Y_1)\otimes V(-Y_2).
      \end{align*}

\item \textrm{Gluing axiom:}
Let $Y_{\mathrm{gl}}$ be the surface obtained from
gluing two boundary components of a surface $Y$.
Then
\[V(Y_{\mathrm{gl}})\cong\bigoplus_{X_i\in A_\mcC} V(Y;(X_i,X_i^*)).\]
 The isomorphism is associative and
compatible with mapping class group projective actions.

Moreover, the isomorphism is compatible with duality as follows:
Let
\begin{align*}
\bigoplus_{j\in \Pi_{\mcC}}\alpha_j &\in V(Y_{\mathrm{gl}};X_{\ell})\cong \bigoplus_{j\in \Pi_{\mcC}}V(Y;X_{\ell},(X_j,X_j^*)),\\
\bigoplus_{j\in \Pi_{\mcC}}\beta_j &\in V(-Y_{\mathrm{gl}};X_{\ell}^*)\cong \bigoplus_{j\in \Pi_{\mcC}} V(-Y;X_{\ell}^*,(X_j,X_j^*)).
\end{align*}
Then there is a nonzero real number $s_{j}$ for each label $j$ such that
\[\biggl\langle \bigoplus_{j\in \Pi_{\mcC}}\alpha_j, \bigoplus_{j\in \Pi_{\mcC}}\beta_j\biggr\rangle =\sum_{j\in \Pi_{\mcC}}
s_j\langle \alpha_j,\beta_j\rangle .\]

\item \textrm{Unitarity:}
Each $V(Y)$ is endowed with a positive-definite Hermitian pairing
\[(\;|\;)  : \overline{V(Y)}\times V(Y)\to \mathbb{C},\]
and each morphism is unitary.  The Hermitian structures are
required to satisfy compatibility conditions as in the duality
axiom. In particular,
\[\biggl( \bigoplus_{i}v_i \biggl| \bigoplus_j w_j\biggr) =\sum_i s_{i0}( v_i|w_i) ,\] for some positive real number $s_{i0}$.
 Moreover, the
following diagram commutes for all $Y$:
$$\xymatrix{
  V(Y) \ar[r]^\cong \ar[d]_\cong & V(-Y)^* \ar[d]^\cong\\
  \overline{V(Y)^*} \ar[r]_\cong & \overline{V(-Y)}}$$
  
\item \textrm{Algebraic axiom:}
There is a choice of basis for all representations with respect to which all matrix entries are algebraic numbers.

\end{enumerate}

\end{definition}

All usual UTMFs satisfy the algebraic axiom \cite{DHW13}.

\subsubsection{Modular categories and their basic invariants}\label{sssection: Modular}

%\subsubsection{Braidings}\label{Braiding}

Given a fusion category $\mcC$, a {\it{braiding}} $c$ of $\mcC$ is a natural family of isomorphisms $c_{V, W}: V \otimes W \to
    W \otimes V$ 
    %in $V$ and $W$ of $\mcC$
    which satisfy the hexagon axioms 
    %for all $U, V, W \in \mcC$ 
    (see \cite{JS93}).  A {\textit{braided fusion category}} is a pair $(\mcC, c)$ in which $c$ is a braiding
    of the fusion category $\mcC$. 
    
For example, the category $\Rep(G)$ of complex representations of a finite group $G$ forms a braided fusion category with the usual tensor product and symmetric braiding $\sigma_{V,W}(v\otimes w)=w\otimes v$.  Any braided fusion category that is equivalent to $\Rep(G)$ for some $G$ as braided fusion categories is called \emph{Tannakian}.  More generally, a braided fusion category with the property that $c_{Y,X}c_{X,Y}=\Id_{X\otimes Y}$ for all $X,Y$ is called \emph{symmetric}.  Non-Tannakian braided fusion categories exist, for example the category $\sVec$ of super-vector spaces with braiding defined on homogeneous objects $V_1,V_2$ by $c_{V_1,V_2}=(-1)^{|V_1|\cdot |V_2|}\sigma_{V_1,V_2}$ where $|V_i|\in\{0,1\}$ is the parity of $V_i$.  An object $X$ for which $c_{Y,X}c_{X,Y}=\Id_{X\otimes Y}$ for all $Y$ is called \emph{transparent}, and the transparent objects generate a symmetric fusion subcategory $\mcC^\prime$ known as the \emph{M\"uger center}.

%-----------------------------------------------------------------------------
% Spherical Fusion Categories
  %-----------------------------------------------------------------------------
%\subsubsection{Spherical Fusion Categories}\label{subsubsection: Spherical Fusion Categories}

Rigidity for a fusion category is encoded in a duality functor $X\rightarrow X^*, f\rightarrow f^*$.  A {\it{pivotal structure}} on a fusion category $\mcC$ is an isomorphism
$j:\id_\mcC \to (-)\bidu$ of monoidal functors, which can then be used to define left and right pivotal traces on endomorphisms $g: V \rightarrow V$ as in \cite[Section 2.1.3]{BNRW151}. 
A pivotal structure is called \emph{spherical} if the left and right pivotal traces coincide, in which case we denote it by $\tr:\End(V)\rightarrow \End(\unit)$.  As $\unit$ is simple we may identify $\tr(f)\in\End(V)$ with a scalar in $\mbbC$.  The \emph{quantum or categorical dimension} of an object $V$ is defined as (the scalar) $d(V):=\tr(\Id_V)$.

%-----------------------------------------------------------------------------
% Modular Categories
%-----------------------------------------------------------------------------
%\subsubsection{Modular Categories}\label{subsubsection: Modular Categories}
A {\it twist}  (or {\it ribbon structure}) of a braided fusion category $(\mcC, c)$ is an $\mbbC$-linear automorphism,
$\theta$, of $\Id_\mcC$ which satisfies
\begin{equation*}
      \theta_{V \otimes W} = (\theta_V \otimes \theta_W)\circ c_{W, V}\circ  c_{V, W}, \quad \theta^{*}_V = \theta_{V^{*}}
\end{equation*}

for $V, W \in \mcC$. A braided fusion category equipped with a ribbon structure is called a {\it ribbon fusion} or {\it premodular} category.
 A premodular category $\mcC$ is called a {\it modular category} if the $S$-{matrix} of $\mcC$, defined by
 \begin{equation*}
      S_{ij} = \tr(c_{V_j, V_{i^*}} \circ c_{V_{i^*}, V_j}) \text{ for }  i, j \in \Pi_\mcC\,,
 \end{equation*}

is non-singular.  Clearly $S$ is a symmetric matrix and $d(V_i)=S_{0i}=S_{i0}$ for all $i$. Note that for a modular category $\mcC$ the M\"uger center $\mcC^\prime\cong\vect$ is trivial, so that we could define modular categories as premodular with $\mcC^\prime\cong\vect$.  Many examples of modular categories will appear in Section \ref{section: Modular Categories}.  

  %-----------------------------------------------------------------------------
  % Grothendieck Ring and Dimensions
  %-----------------------------------------------------------------------------
\subsubsection{Grothendieck ring and dimensions}\label{groth and dim}
The {\it Grothendieck ring} $K_0(\mcC)$ of a fusion category $\mcC$ is the based $\mbbZ$-ring \cite{OstrModule} generated by $\Pi_\mcC$ with multiplication induced from $\otimes$.  The structure coefficients of $K_0(\mcC)$ are obtained from:
\begin{equation*}
      V_i \otimes V_j \cong \bigoplus_{k \in \Pi_\mcC} N_{i,j}^k \,V_k
\end{equation*}

where $N_{i,j}^k = \dim(\Hom_{\mcC}(V_k, V_i\otimes V_j))$. This family of non-negative
integers $\{N_{i,j}^k\}_{i,j,k \in \Pi_\mcC}$ is called the \textit{fusion
    rules} of $\mcC$.
Two fusion categories $\mcC$ and $\mcD$ are called \emph{Grothendieck equivalent} if $K_0(\mcC)\cong K_0(\mcD)$ as $\mbbZ$-rings.  
For a braided fusion category $K_0(\mcC)$ is a commutative ring and the fusion rules satisfy the symmetries:
\begin{equation}\label{fusion symmetries}
      N_{i,j}^k=N_{j,i}^k=N_{i,k^*}^{j^*}=N_{i^*,j^*}^{k^*},\quad
      N_{i,j}^0=\delta_{i,j^*}
\end{equation}

The {\it fusion matrix} $N_i$ of $V_i$, defined by $(N_i)_{k,j}=N_{i,j}^k$, is an integral matrix with non-negative entries. In the braided fusion setting these matrices are normal and mutually commuting. The largest real eigenvalue of $N_i$ is called the {\it Frobenius-Perron dimension} of $V_i$ and is denoted by $\FPdim(V_i)$.
Moreover, $\FPdim$ can be extended to a $\mbbZ$-ring homomorphism from $K_0(\mcC)$ to $\mbbR$ and is the unique such homomorphism that is positive (real-valued) on $\Pi_\mcC$ (see \cite{ENO1}). The {\it Frobenius-Perron dimension} of $\mcC$ is defined as
\begin{equation*}
      \FPdim(\mcC) = \sum_{i \in \Pi_\mcC} \FPdim(V_i)^2\,.
\end{equation*}

\begin{definition}
 A fusion category $\mcC$ is said to be
 \begin{enumerate}
 \item {\it weakly integral} if $\FPdim(\mcC)\in\mbbZ$.
 \item {\it integral} if $\FPdim(V_{j})\in\mbbZ$ for all $j\in\Pi_{\mcC}$.
 \item {\it pointed} if $\FPdim(V_{j})=1$ for all $j\in\Pi_{\mcC}$.
 \end{enumerate}

Furthermore, if $\FPdim(V)=1$, then $V$ is {\it invertible}.
\end{definition}
The \emph{categorical dimension} $\dim(\mcC)$ of $\mcC$ is defined similarly: $$\dim(\mcC)=\sum_{i\in\Pi_\mcC}d(V_i)^2$$ and if $\FPdim(\mcC)=\dim(\mcC)$ we say that $\mcC$ is \emph{pseudo-unitary}.

The columns of the $S$-matrix can be seen to be simultaneous eigenvectors for the (commuting) fusion matrices $\{N_a:a\in \Pi_\mcC\}$, so that $S$ diagonalizes all $N_a$ simultaneously.  This leads to the famous \textit{Verlinde formula}: $N_{ij}^k=\frac{1}{\dim(\mcC)}\sum_r\frac{S_{ir}S_{jr}S_{\hat{k}r}}{S_{0r}}$.

\subsection{From Modular Categories to Modular Functors and Back}

TQFTs such as the Reshetikhin-Turaev ones can have framing anomaly which are natural for applications to physics for chiral phases with non-zero central charges.  We follow the physical solution to consider only quantum states projectively, i.e. rays in Hilbert spaces of states.  For more detailed introduction to TQFTs, see \cite{Tu94,Freed13}.

\subsubsection{(2+1)-TQFTs}

The oriented bordism category $\textrm{Bord}^3$ has closed oriented surfaces as objects and equivalence classes of compact oriented bordisms between them as morphisms.

We enlarge the projective category $\Hil$ of finite dimensional Hilbert spaces in section \ref{Sec:UTMF} to include partial unitaries as morphisms.

\begin{definition}

A unitary $(2+1)$-TQFT with a strict unitary fusion category $\mcC$ is a pair $(Z,V)$, where $Z$ is a unitary symmetric monoidal projective functor from $\textrm{Bord}^3$ to $\Hil$, and $V$ is a UTMF with $\mcC$ from $\textrm{Bord}_{(2,\mcC)}$ to $\Hil$.  We will call $Z$ the partition functor.  The functors $Z$ and $V$ are compatible in the sense that when oriented closed surfaces and their mapping classes are included into $\text{Bord}^3$ via the mapping cylinder construction, they are functorially equivalent.

\end{definition}

\begin{theorem}

\begin{enumerate}
\item Given a modular category $\mcB$, there is a $(2+1)$-TQFT $(Z,V)$ and $V$ is based on the fusion category underlying $\mcB$.

\item Given a modular functor $V$ associated with a fusion category $\mcC$, there is a $(2+1)$-TQFT $(Z,V,\mcB)$ with modular category $\mcB$ extending $\mcC$ in the sense that $\mcB$ is equivalent to $\mcC$ as fusion categories.
\end{enumerate}

\end{theorem}

This theorem is from \cite{Tu94}.
The subtlety for the correspondence about TQFTs and modular functors disappears because our input data for a modular functor is a strict fusion category, instead of a monoidal class.

\subsubsection{Dictionary of modular categories and anyon systems}

A dictionary of terminologies between UMCs and anyons systems is given in Table~\ref{table-dictionary} \cite{Wa10}.  The triangular space $V_{ab}^c$ is $\textrm{Hom}(a\otimes b,c)$ and $V(Y)$ the Hilbert space that the TQFT assigned to the surface $Y$.  Tangles appear because of the creation and annihilation of pairs of anyons.  For $F$-matrices and the relation between anyons and quantum topology, we refer to Chapter 6 of \cite{Wa10} for more detail.
\begin{table}[hb]
\begin{tabular}{|l|l|}\hline
\emph{UMC} & \emph{Anyonic system}\\\hline
simple object & anyon\\\hline
label & anyon type or topological charge\\\hline
tensor product & fusion\\\hline
fusion rules & fusion rules\\\hline
triangular space $V^c_{ab}$ or $V^{ab}_c$ & fusion/splitting space\\\hline
dual & antiparticle\\\hline
birth/death & creation/annihilation\\\hline
mapping class group representations & anyon statistics\\\hline
nonzero vector in $V(Y)$ & ground state vector\\\hline
unitary $F$-matrices & recoupling rules\\\hline
twist $\theta_x = e^{2\pi i s_x}$ & topological spin\\\hline
morphism & physical process or operator\\\hline
tangles & anyon trajectories\\\hline
quantum invariants & topological amplitudes\\ \hline
\end{tabular}
\caption{}
\label{table-dictionary}
\end{table}

\subsection{Graphical Calculus}

A powerful tool to study tensor categories with structures is the graphical calculus---a far reaching generalization of spin networks.
It is important to understand the physical interpretation of graphical calculus using anyon theory.  Basically all categorical axioms and structures in UMCs make the graphical calculus work out as expected when care is taken.

\subsubsection{Handle and sign of an anyon}

\begin{figure}[h]\includegraphics[width=2.7in]{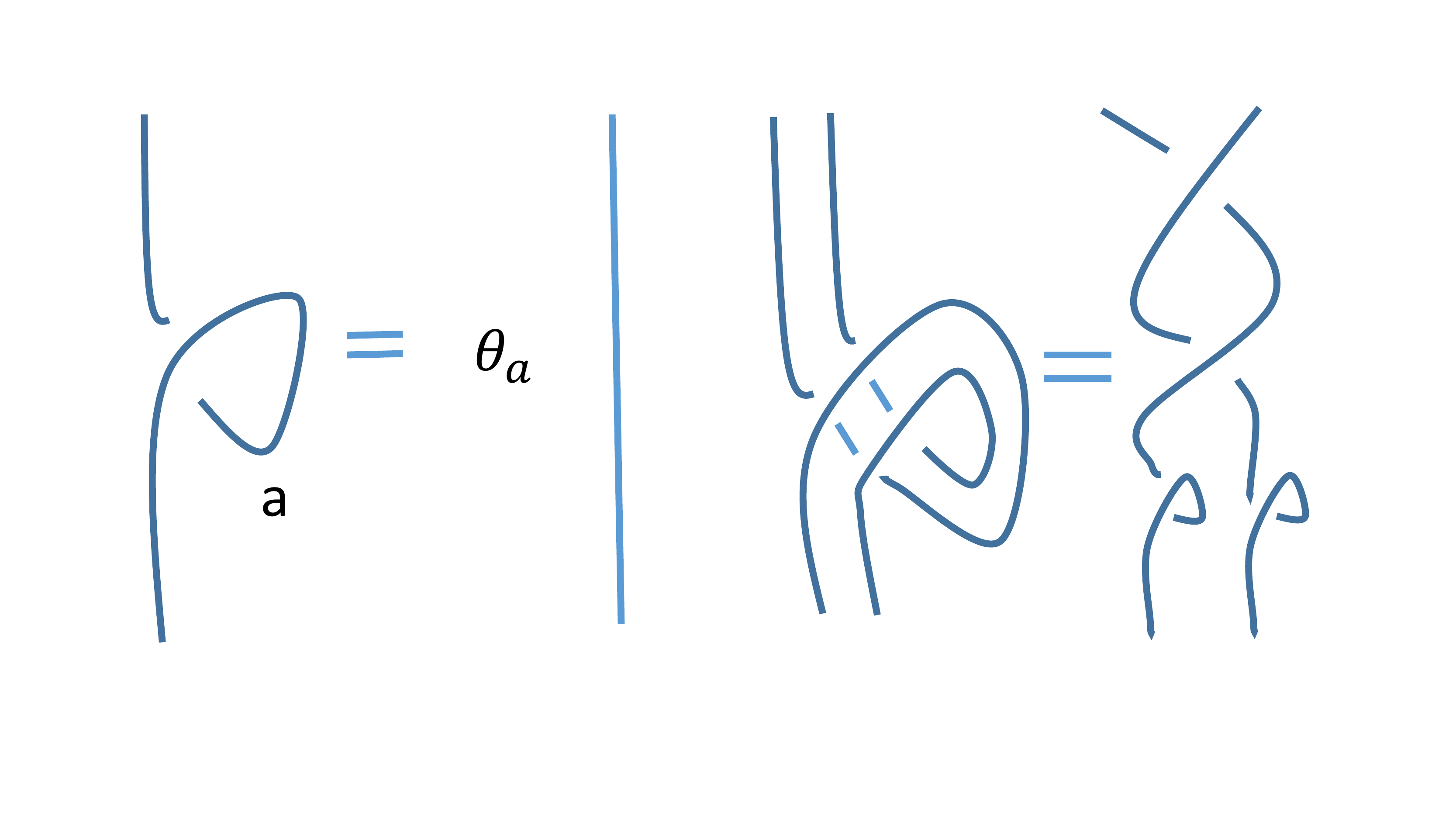}\hspace{1pc}
\begin{minipage}[b]{20pc} \caption{\label{spinstats} Spin-Statistics Connection: Here $\theta_a=e^{2\pi ih_a}$ for some $h_a\in\mbbQ$ is the topological twist.  The picture on the left acquires the topological twist $\theta_a$ when pulled straight, not as in knot diagrams. }\vspace{.1in}
\end{minipage}
\end{figure}

Classical particles are perceived as points, therefore their configuration spaces are generally manifolds.  A defining feature of a quantum particle such as the electron is its spin.  A deep principle is the spin-statistics connection.  It is common to picture the spin of a particle through a $2\pi$ rotation as in Figure \ref{spinstats}.  The $h_a$'s mod $1$ in the topological spin correspond to scaling dimensions if there is a corresponding conformal field theory.  The existence of topological spin means that knot diagrams are equivalent only up to framed link diagrams.  The right-handed picture represents the compatibility of tensor product, twist and braidings.

To encode the topological twist in graphical calculus, we use ribbons to representation anyon trajectories, i.e. we consider framed knot diagrams.  Pulling the ribbon tight results in a full twist of one side around the other side of the ribbon.  The picture of a quantum particle would be a small arrow\footnote{One justification is as follows.  A quantum particle is an elementary excitation in a quantum system, so is represented by a non-zero vector in a Hilbert space.  Even after the state vector is normalized, there is still a phase ambiguity $e^{i\theta}, \theta \in [0,2\pi)$, which parameterizes the standard circle in the complex plane.  So a semi-classical picture of the particle sitting at the origin of the plane would be an arrow from the origin to a certain angle $\theta$ (this arrow is really in the tangent space of the origin), which also explains the base-point on the boundary of the infinitesimal disk.}.  Visualizing a quantum particle as an infinitesimal arrow leads to the ribbon picture of the worldline of the quantum particle.

An anyon has spin so it needs a handle, which requires a base-point for the boundary circles.  The handle has no categorical interpretation.  We need a sign to indicate if the anyon comes from the past or the future, and it corresponds to the dual objects.  Therefore, a signed object in a strict fusion category is used to color strings in the graphical calculus. A complete treatment of this diagrammatic yoga of \textit{graphical calculus} involving such pictures can be found in \cite[Appendix E]{Ki06},\cite{Hagge08}, and \cite[Section 4.2]{Wa10}.

\subsection{Anyonic Quantum Computing Models}

\begin{figure}[ht]\includegraphics[width=5.7in]{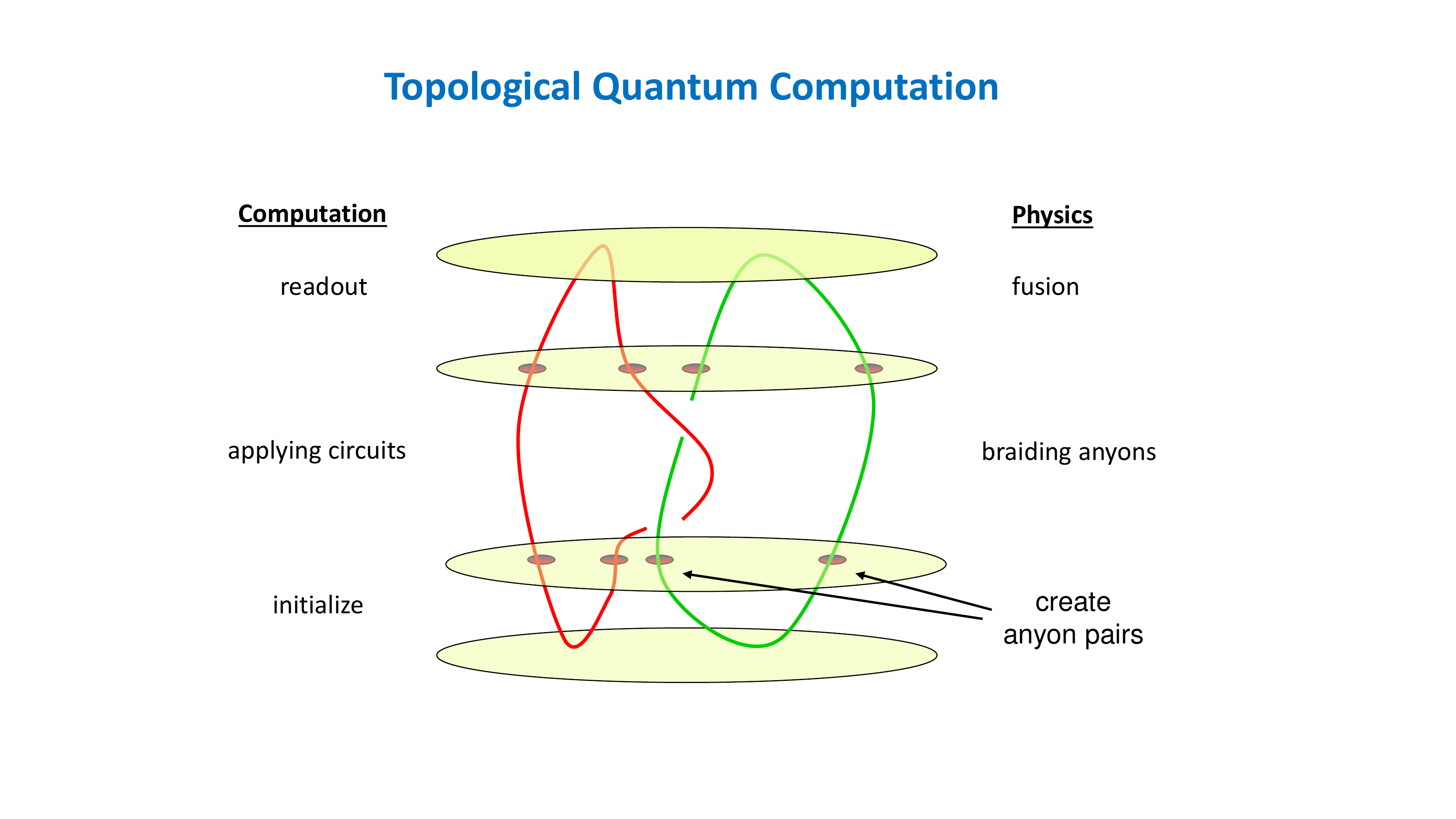}
 \caption{\label{fig:TQC} TQC Setup: time flows from bottom to top.}

\end{figure}

This basic setup of TQC is depicted in Fig.~\ref{fig:TQC}.

So anyonic computing runs as follows: 1) one starts the computation by creating anyon pairs from the vacuum to encode the input, 2) braiding these anyons to move the initial state, and 3) measuring the anyon type of pairs of neighboring particles.  Categorically, the first step is the implementation of a morphism in $\textrm{Hom}(\unit, X^{\otimes n})$ for some simple object $X$ and $n$.  The second step is a braiding.  The third measurement step is an implementation of a morphism in $\textrm{Hom}(X^{\otimes n},\unit)$. One key is that after braiding the anyons, a neighboring anyon pair may have obtained a different total charge (besides $0$, i.e. the vacuum), i.e. a simple object in the fusion result.  The computing result is a probability distribution on anyon types obtained by repeating the same process polynomially many times, taking a tally of the output anyon types.  Usually the probability of neighboring pairs of anyons returning to vacuum is taken as the computing answer, which approximates some topological invariant of links obtained from the braiding process (the links are the plat closures of the braids).  

The topological invariance ensures that slight variations in the process (e.g. small deviations in the trajectory of an anyon in space-time) do not influence the output.
 
The time evolution of the topological state spaces $V(Y;\tau_i)$ must be a unitary operator.  In particular, a sequence of anyon braidings corresponding to a braid $\beta$ induces a unitary transformation
$|\psi\rangle\mapsto U_\beta|\psi\rangle$.  In the topological model of \cite{FKLW03} these are the \textit{quantum circuits}. 

Informally, an anyon $X$ is called \textit{braiding universal} if any computation can be approximately achieved by braiding multiple $X$'s (more detail can be found in section \ref{Sec:AQC}).  Since universality is of paramount importance of quantum computation, it is very important to know if a particular anyon is braiding universal.

\subsection{Anyons In The Real World}

TPMs in nature include quantum Hall states--integral and fractional---and the recently discovered topological insulators \cite{NSSFD08,Majorana15}.  One incarnation of the Ising anyon $\sigma$ is the Majorana zero mode, which is experimentally pursued in nanowires \cite{Mou12,Marcus16}.

\begin{figure}[ht]
\includegraphics[width=6in]{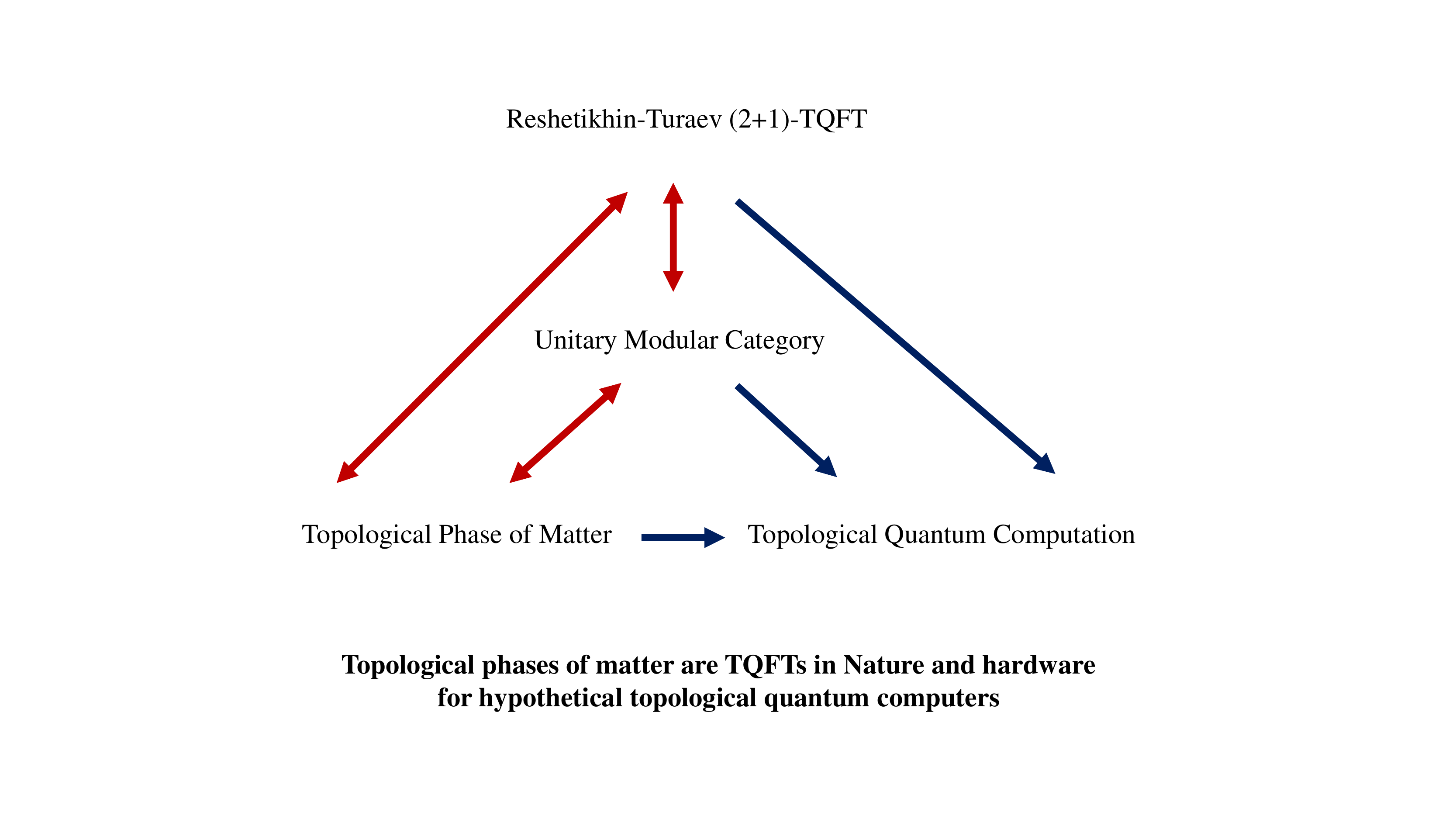}
\caption{\label{figNexus} The Nexus.}
\end{figure}

The nexus among TQC, TPM and TQFT is summarized in  Figure \ref{figNexus}.  Topological phases of matter are states of matter whose low energy physics are modeled by unitary TQFTs, and they can be used to construct large scale topological quantum computers.

\section{Computational Power of Physical Theories}

Constructing machines is a defining characteristic of human.  Every new physical theory provides an opportunity to build new kinds of computing machines.  Freedman articulated this idea in \cite{Fr98}: \lq\lq As a generality, we propose that each physical theory supports computational models whose power is limited by the physical theory.  It is well known that classical physics supports a multitude of the implementation of the Turing machine".
Freedman further suggested that computational models based on some TQFTs might be more powerful than quantum computing---the computing model based on quantum mechanics.  But when accuracy and measurement are carefully analyzed, the computing model based on TQFTs are polynomially equivalent to quantum computing \cite{FKLW03}.  This mathematical theorem  suggests the possibility that any computing model implemented within quantum field theory is polynomially equivalent to quantum computing.  But we do not have mathematical formulations of quantum field theories to prove that hypercomputation is impossible.

In this section we give an introduction to von Neumann's axiomatization of quantum mechanics and quantum computing.  Abstract quantum mechanics for finite dimensional systems is completely elementary, in striking contrast to quantum field theory.  More or less, it is simply a physical way of thinking about complex linear algebra.

\subsection{Encoding and Computing Problems}

\begin{definition}

A computing problem is a sequence of Boolean functions $f_n: \mbbZ_2^n \rightarrow \mbbZ_2^{a(n)}$.
A computing problem with $a(n)=2$ is called a decision problem.
A reversible computing problem is one such that $a(n)=n$ and every $f_n$ is a permutation.

Let $\{0,1\}^{*}=\cup_{n=0}^{\infty}\mbbZ_2^n$, where for $n=0$, $\mbbZ_2^n$ is the empty string.  Then a computing problem is simply a map $f: \{0,1\}^{*} \rightarrow \{0,1\}^{*}$. 

\end{definition}

To formalize concrete computing problems such as  integer factoring or Jones polynomial evaluation into families of Boolean maps, it is important that we encode the inputs such as integers and links into bit strings in an intelligent way because encodings can change computational complexities.  If an integer $N$ is encoded by unary strings, then dividing $N$ by primes is polynomial in the input length.  Most other reasonable encodings will lead to the same complexity class, therefore, we will not discuss encodings further (see \cite{garey02}).

We consider classical physics as part of quantum physics, and will embed classical computation into quantum computation through reversible classical computing.  Therefore, quantum computing can solve the same class of computing problems potentially much faster.  Note that quantum computers can be simulated by classical computers, though potentially exponentially slowly.

\subsection{Quantum Framework}

Some basic physical principles for quantum systems include: there is a definite state at each moment, and its evolution from one moment to another is deterministic.  Moreover, a known state can be prepared and a known evolution can be repeated as many times as we need.

Quantum mechanics is a set of rules that predict the responses of the microscopic world to our measuring devices.  The most salient feature is the superposition of different states.  With the advent of quantum information science,  another quantum correlation comes to the center stage: entanglement--{\it the} characteristic attribute of quantum mechanics according to Schr\"odinger.

von Neumann's axiomatization of quantum theory consists of four principles: superposition, linear evolution, entanglement, and projective measurement.  The controversial\footnote{The interpretation of measurement in quantum mechanics is still under debate, and the projective measurement here is the one usually used in quantum computation.} measurement reflects well the un-controllable disturbance of the quantum state by our measuring devices: energy injecting into the quantum system by measuring apparatus overwhelms the fragile state.

Two operations on Hilbert spaces that are used to describe superposition and entanglement, respectively, are the direct sum $\oplus$ and the tensor product $\otimes$. While interference is arguably more fundamental, a deeper understanding of quantum mechanics would come from the interplay of the two.  The role for entanglement is more pronounced for many-body quantum systems such as systems of $10^{11}$ electrons in condensed matter physics.  Advancing our understanding of the role of entanglement for both quantum computing and condensed matter physics lies at the frontier of current research.

\subsubsection{von Neumann axioms}
von Neumann's axioms for quantum mechanics are as follows:
\begin{enumerate}
\item \emph{State space:}
There is a Hilbert space $L$ describing all possible states of a quantum system.
Any nonzero vector $|v\rg$ represents a state, and two nonzero vectors $|v_1\rg$ and $|v_2\rg$ represent the same state iff $|v_1\rg = \lambda |v_2\rg$ for some scalar $\lambda \neq 0$.
Hilbert space embodies the superposition principle.

Quantum computation uses ordinary finite-dimensional Hilbert space $\mbbC^m$, whose states correspond to projective points ${\mbbC P}^{m-1}$.
Therefore information is stored in state vectors, or more precisely, points on ${\mbbC P}^{m-1}$.

\item \emph{Evolution:} If a quantum system is governed by a quantum Hamiltonian operator $H$, then its state vector $|\psi\rg$  is evolved by solving the  Schr\"odinger equation $i\hbar \tfrac{\partial |\psi\rg}{\partial t} = H|\psi\rg$.
When the state space is finite-dimensional, the solution is $|\psi_t\rg = e^{-\frac{i}{\hbar}tH}|\psi_0\rg$ for some initial state $\ket{\psi_0}$.
Since $H$ is Hermitian, $e^{-\frac{i}{\hbar}tH}$ is a unitary transformation.
Therefore we will just say states evolve by unitary transformations.

In quantum computation, we apply unitary transformations to state vectors $\ket{\psi}$ to process the information encoded in $\ket{\psi}$.
Hence information processing in quantum computation is multiplication by unitary matrices.
\item \emph{Measurement:} Measurement of a quantum system is given by a Hermitian operator $M$ such as the Hamiltonian (= total energy).
Since $M$ is Hermitian, its eigenvalues are real.
If  they are pairwise distinct, we say the measurement is complete.
Given a complete measurement $M$ with eigenvalues $\{\lambda_i\}$, let $\{e_i\}$ be an orthonormal basis of eigenvectors of $M$ corresponding to $\{\lambda_i\}$.
If we measure $M$ in a normalized state $|\psi\rg$, which can be written as $|\psi\rg = \sum_ia_i|e_i\rg$, then the system will be in state $|e_i\rg$ with probability $|a_i|^2$ after the measurement. Conceptually,  quantum mechanics is like a square root of probability theory because amplitudes are square roots of probabilities.
The basis $\{e_i\}$ consists of states that are classical in a sense.
This is called projective measurement and is our read-out for quantum computation.

Measurement interrupts the deterministic unitary evolution and outputs a random variable $X \colon \{e_i\} \to \{\lambda_i\}$ with probability distribution $p(X=\lambda_i) = |a_i|^2$, and hence is the source of the probabilistic nature of quantum computation.

\item \emph{Composite system:} If two systems with Hilbert spaces $L_1$ and $L_2$ are brought together, then the state space of the joint system is $L_1 \otimes L_2$.
\end{enumerate}

\subsubsection{Entanglement}

Superposition is only meaningful if we have a preferred basis, i.e. a direct sum decomposition of the Hilbert space.  
Entanglement is a property of a composite quantum system, so is only meaningful with a tensor product decomposition.  

\begin{definition}

Given a composite quantum system with Hilbert space $\mcL=\bigotimes_{i\in I} L_i$, then a state is entangled with respect to this tensor decomposition into subsystems if it is not of the form $\otimes_{i\in I}|\psi_i\rangle$ for some $|\psi_i\rangle \in L_i$.

\end{definition}

Quantum states that are not entangled are product states, which behave as classical states.  Entangled states are composite states for which the constituent subsystems do not have definite states and serve as sources of quantum weirdness---Einstein's spooky action at distance.

\subsubsection{Error-correcting codes}

Recall that an operator $\mcO$ is $k$-local for some integer $k\geq 0$ on some qubits $(\mbbC^2)^{\otimes m}, m>k,$ if $\mcO$ is of the form $\textrm{Id}\otimes A  \otimes \textrm{Id}$, where $A$ acts on $k$ qubits.  

An error-correcting code is an embedding of $(\mbbC^2)^{\otimes n}$ into $(\mbbC^2)^{\otimes m}$ such that information in the image of $(\mbbC^2)^{\otimes n}$ is protected from local errors on $(\mbbC^2)^{\otimes m}$, i.e. $k$-local operators for some $k$ on $(\mbbC^2)^{\otimes m}$ cannot change the embeded states of $(\mbbC^2)^{\otimes n}$ in $(\mbbC^2)^{\otimes m}$, called the code subspace, in an irreversible way.
We call the encoded qubits the logical qubits and the raw qubits $(\mbbC^2)^{\otimes n}$ the physical qubits.

The following theorem can be found in \cite{Gott00}, which can also be used as a definition of error correcting code:

\begin{theorem}
Let $V,W$ be logical and physical qubit spaces.  
The pair $(V,W)$ is an error-correcting code if there exists an integer $k \geq 0$ such that the composition
\[ \xymatrix{ V \ar@{^(->}[r]^i & W \ar[r]^{O_k} & W \ar@{->>}[r]^\pi & V }\]
is $\lambda \cdot \id_V$ for any $k$-local operator $O_k$ on $W$, where $i$ is inclusion and $\pi$ projection.
\end{theorem}

When $\lambda \neq 0$, $O_k$ does not degrade the logical qubits as projectively its action on the logical qubits is the identity operator.
But when $\lambda=0$,  $O_k$ rotates logical qubits out of the code subspace, introducing errors.
But such local errors always rotate a state in the code subspace $V$ to an orthogonal state $W$, it follows that local errors can be detected by simultaneous measurements and subsequently corrected by using designed gates.  It follows that local errors are detectable and correctable.

The possibility of fault-tolerant quantum computation was a milestone in quantum computing.
The smallest number of physical qubits fully protecting one logical qubit is $5$.

\subsection{Gates, Circuits, and Universality}

A gate set $\mc{S}$ is the elementary operations that we will carry out repeatedly to complete a computational task.
Each application of a gate is considered a single step, hence the number of gate applications in an algorithm represents consumed time, and is a complexity measure.
A gate set should be physically realizable and complicated enough to perform any computation given enough time.
It is not mathematically possible to define when a gate set is physical as ultimately the answer comes from physical realization.
Considering this physical constraint, we will require that all entries of gate matrices are efficiently computable numbers when we define complexity classes depending on a gate set.
Generally, a gate set $\mc{S}$ is any collection of unitary matrices in $\bigcup_{n=1}^\infty\mathrm{U}(2^n)$.
Our choice is
\[ \mc{S} = \{H, \sigma_z^{1/4}, \mathrm{CNOT} \} \]
where
\begin{align*}
H &= \frac{1}{\sqrt{2}} \begin{pmatrix}
      1 & 1\\
      1 & -1
    \end{pmatrix} &\text{is the Hadamard matrix},\\
\sigma_z^{1/4} &= \begin{pmatrix}
      1 & 0\\
      0 & e^{\pi i/4}
    \end{pmatrix} &\text{is called the $\frac{\pi}{8}$ or T-gate},\\
\mathrm{CNOT} &= \begin{pmatrix}
    1&0&0&0\\
    0&1&0&0\\
    0&0&0&1\\
    0&0&1&0
  \end{pmatrix} &\text{ in the two-qubit basis $\bigl\{\ket{00},\ket{01},\ket{10},\ket{11}\bigr\}$.}
\end{align*}
The CNOT gate is called controlled-NOT because the first qubit is the control bit, so that when it is $\ket{0}$, nothing is done to the second qubit, but when it is $\ket{1}$, the NOT gate is applied to the second qubit.

\begin{definition}\
\begin{enumerate}
\item
An $n$-qubit quantum circuit over a gate set $\mc{S}$ is a map $U_L \colon (\mbbC^2)^{\otimes n} \to (\mbbC^2)^{\otimes n}$ composed of finitely many matrices of the form $\idmat_p \otimes g \otimes \idmat_q$, where $g \in \mc{S}$ and $p,q$ can be $0$.
\item A gate set is universal if the collection of all $n$-qubit circuits forms a dense subset of $\SU(2^n)$ for any $n$.
\end{enumerate}
\end{definition}
The gate set $\mc{S} = \{H, \sigma_z^{1/4}, \mathrm{CNOT}\}$ will be called the standard gate set, which we will use unless stated otherwise.

\begin{theorem}\label{thm-univgates}\
\begin{enumerate}
\item The standard gate set is universal.
\item Every matrix in $\mbbU(2^n)$ can be efficiently approximated up to an overall phase by a circuit over $\mc{S}$.
\end{enumerate}
\end{theorem}

Note that (2) means that \emph{efficient} approximations of unitary matrices follows from approximations for free due to the Solvay-Kitaev theorem. See \cite{NC00}.

\subsection{Complexity Class BQP}

Let $\mcC_T$ be the class of computable functions $\mcC_T$ and $P$ the subclass of functions efficiently computable by a Turing machine\footnote{More precisely, this should be denoted as FP as P usually denotes the subset of decision problems.}.
Defining a computing model $\mcX$ is the same as selecting a class of computable functions from $\mcC_T$, denoted as $\mcX P$.  The class $\mcX P$ of efficiently computable functions codifies the computational power of computing machines in $\mcX$.  Quantum computing selects a new class BQP--bounded-error quantum polynomial time. The class BQP consists of those problems that can be solved efficiently by a quantum computer. In theoretical computer science, separation of complexity classes is extremely hard as the millennium problem $P$ vs $NP$ problem shows.  It is generally believed that the class BQP does not contain NP-complete problems.  Therefore, good target problems will be those NP problems which are not known to be NP complete.  Three candidates are factoring integers, graph isomorphism, and finding the shortest vector in lattices.

\begin{definition}
Let $\mathcal{S}$ be any finite universal gate set with efficiently computable matrix entries.
A problem $f : \{0,1\}^* \to \{0,1\}^*$ \textup{(}represented by $\{f_n\}: \mbbZ_2^n \to \mbbZ_2^{m(n)}$\textup{)}
is in BQP \textup{(}i.e., can be solved efficiently by a quantum computer\textup{)}
if there exist polynomials $a(n),g(n) : \mbbN \to \mbbN$ satisfying $n+a(n) = m(n)+g(n)$
and a classical efficient algorithm to output a map $\delta(n) : \mbbN \to \{0,1\}^*$ describing a quantum circuit $U_{\delta(n)}$ over $\mathcal{S}$ of size $O(\mathrm{poly}(n))$ such that:
\begin{align*}
U_{\delta(n)}\ket{x,0^{a(n)}} &= \sum_I a_I\ket{I}\\
\sum_{\ket{I} = \ket{f_n(x),z}}\!\!\!\! |a_I|^2 &\geq \frac{3}{4}, \quad \text{where $z \in \mbbZ_2^{g(n)}$}
\end{align*}
\end{definition}
The $a(n)$ qubits are an ancillary working space, so we initialize an input $\ket{x}$ by appending $a(n)$ zeros and
identify the resulting bit string as a basis vector in $(\mbbC^2)^{\otimes (n+a(n))}$.
The $g(n)$ qubits are garbage.
The classical algorithm takes as input the length $n$ and returns a description of the quantum circuit $U_{\delta(n)}$.
For a given $\ket{x}$, the probability that the first $m(n)$ bits of the output equal $f_n(x)$ is $\geq \frac{3}{4}$.

The class BQP is independent of the choice of gate set as long as the gate set is efficiently computable.
The threshold $\frac{3}{4}$ can be replaced by any constant in $(\frac{1}{2},1]$.

Classical computation moves the input $x$ through sequences of bit strings $x'$ before we reach the answer bit string $f(x)$.  In quantum computing, the input bit string $x$ is represented as a basis quantum state $| x\rangle\in {(\mbbC^2)}^{\otimes n}$---the Hilbert space of $n$-qubit states.  Then the computing is carried out by multiplying the initial state vector $|x\rangle$ by unitary matrices from solving the Schr\"odinger equation.  The intermediate steps are still deterministic, but now through quantum states which are superpositions of basis states.  The answer will be contained in the final state vector of the quantum system---a superposition with exponentially many terms.  But there is an asymmetry between the input and output because in order to find out the answer, we need to measure the final quantum state.  This quantum measurement leads to a probability by the von Neumann formulation of measurement.

Shor's astonishing algorithm proved that factoring integers is in BQP, which launched quantum computing.  It follows from the above discussion that classical reversible computation is a special case of quantum computing.  Therefore, $\textrm{P}\subseteq \textrm{BQP}$.  Quantum states are notoriously fragile: the fragility of qubits--decoherence--so far has prevented us from building a useful quantum computer.  

\begin{figure}[h]\includegraphics[width=4.1in]{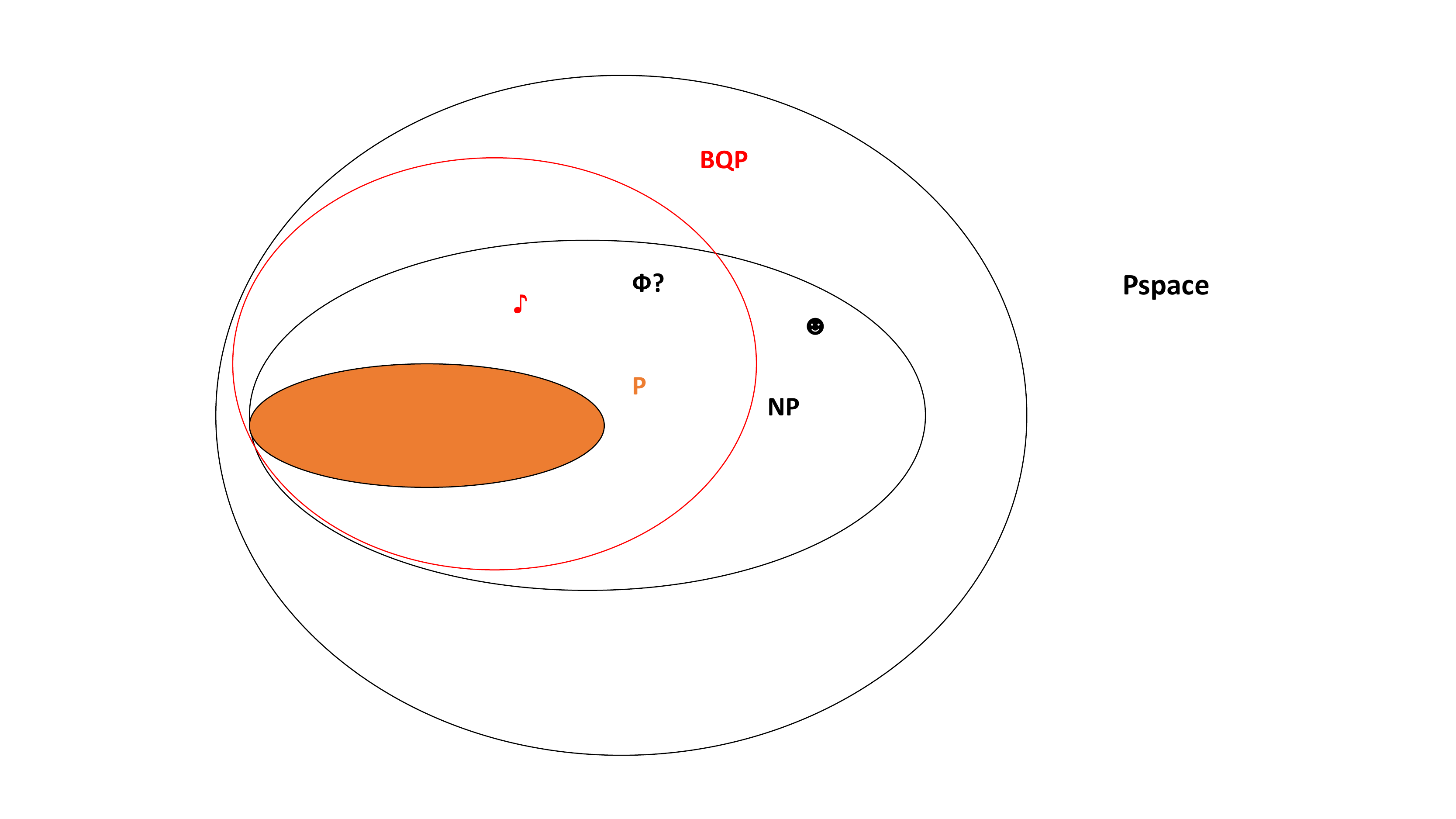}
\begin{minipage}[b]{20pc} \caption{\label{BQP} P vs BQP}
\end{minipage}
\end{figure}

A conjectural picture of the various complexity classes is illustrated in Figure \ref{BQP}, where $\Phi?$ is the graph isomorphism problem, \textmusicalnote\/ is integer factoring and \textbullet\/ represents NP completeness. Mathematically all proper inclusions might collapse to equality.

\subsection{Simulation of TQFTs}\label{Sec:Simulation}

\subsubsection{Hidden locality of TQFT}

The topological Hilbert spaces $V(Y;\tau_i)$ of a TQFT have no natural tensor product structures.  To simulate the mapping class groups action on $V(Y;\tau_i)$ by a quantum computer,  we will embed $V(Y;\tau_i)$ into the tensor product of some qudits $\mathbb{C}^d$ for some integer $d\geq 2$.  The natural choice of qudit is $W=\bigoplus_{(a,b,c\in \Pi_{\mcC})}\Hom(a \otimes b,c)$, then $V(Y;\tau_i)$ embeds into $W^{\otimes m}$ for some $m$ by the gluing formula.  The choice of qudit $W$ is natural because using pant-decomposition of surfaces and the gluing formula, every state in $V(Y;\tau_i)$ is embeded in the tensor products of qudits.

Recall a mapping class is an equivalence class of self-diffeomorphisms of a surface up to isotopy and a Dehn twist is a generating mapping class (for a picture, see \cite{FKW02}).  
To simulate a mapping class unitary $U(f)$ on $V(Y;\tau_i)$, we look for an efficient circuit $U_L$ so that the following diagram commutes:
$$\xymatrix{
V(Y) \ar[d]^{U(f)} \ar[r] & W^{\otimes m} \ar[d]^{U_L}\\
V(Y) \ar[r] &W^{\otimes m} }$$

\begin{theorem}\cite{FKW02}
Using representation matrices of the standard Dehn twists as gates for a quantum computer, we can simulate the mapping class group representation matrices of any unitary $(2+1)$-TQFT efficiently up to a phase.
\end{theorem}

\subsubsection{Approximate topological invariants}

As corollaries of the quantum simulation of TQFTs, we obtain efficient quantum algorithm to approximate link and $3$-manifold invariants.  They are native computing problems for TQC. For example, we have the following from \cite{FKW02,FKLW03}:
\begin{theorem}\label{jonesbqp}
Approximation of the normalized Jones polynomials evaluated at roots of unity $q=e^{\frac{2\pi i}{r}}$ for any $r\geq 3$ is in BQP.
\end{theorem}

\subsection{Open Problems}

Given any mathematically defined physical theory $H$, we can ask if $H$ can be simulated efficiently by quantum computers, and if a new model of computation can be constructed from $H$.  Mathematically formulated quantum field theories are relatively rare, including TQFTs and conformal field theories (CFTs).

\subsubsection{Simulating CFTs}

The quantum polynomial Church-Turing thesis suggests that quantum computers can efficiently simulate CFTs.  But even how to formulate the problem is very difficult \cite{Moji18}.

\section{Topological Phases of Matter}

In this section, we focus on the algebraic and topological study of intrinsic topological order with long-range entanglement in $2D$ topological phases of matter (TPMs).  Mathematically, we are studying UMCs and TQFTs and their realization by lattice models.  Our goal is to lay the mathematical foundations for a mathematical study of TPMs via the Hamiltonian formulation as inspired by condensed matter physics \cite{Wenzoo,Wenbook}.
Microsoft Station Q has been instrumental in many physical developments in TPMs and their experimental realization during the last decade\footnote{
In 2003, Freedman, Nayak, and the second author organized the workshop {\it topology in condensed matter physics} in the American Institute of Mathematics.  Afterwards, Freedman proposed the establishment of a Microsoft research institute on the campus of UC Santa Barbara to build a topological quantum computer.  In 2005, Microsoft Station Q began with Freedman, Kitaev, Nayak, Walker, and the second author.}.  Our main contribution is a mathematical definition of $2D$ TPMs, and the formulation of well-known physical statements into mathematical conjectures.  We omit the analytic aspects of TPMs, symmetry protected topological order and short-range entangled states.

Two guiding principles to study TPMs are locality and unitarity.  In \cite{Fr01}, Freedman wrote:\lq\lq Nature has the habit of intruding on the prodigies of purest thought and encumbering them with unpleasant embellishments.  So it is astonishing when the chthonian hammer of the engineer resonates precisely to the gossamer fluttering of theory.  Such a moment may soon be at hand in the practice and theory of quantum computation.  The most compelling theoretical question, \lq localization,' is yielding an answer which points the way to a solution of Quantum Computing's (QC) most daunting engineering problem: reaching the accuracy threshold for fault tolerant computation."  Our definition of TPMs is a form of localization of TQFTs.  As we saw in Section \ref{Sec:Simulation}, the efficient simulation of TQFTs by quantum computers can be regarded as another localization of TQFTs.  For braid groups, localization has been intensively discussed in \cite{RWlocal}.   While a localization of a TQFT is some progress towards its physical realization, the devil in Nature's \lq\lq unpleasant embellishments" is formidable for the construction of a topological quantum computer.

Hamiltonians in this section are all quantum, so they are Hermitian matrices on Hilbert spaces.  Unless stated otherwise, all our Hilbert spaces are finite dimensional, hence  isomorphic to $\mbbC^d$ for some integer $d$.  We will refer to $\mbbC^d$ as a {\it qudit} following the quantum computing jargon.

\subsection{Quantum Temperature}

Quantum phases of matter are determined by quantum effects, not by the energy {\it alone}. So ideally they are phases at zero temperature.  But in realty, zero temperature is not an option. Then what is a quantum system?  Particles normally have a characteristic oscillation frequency depending on their environment.  By quantum mechanics, they have energy $E=h\nu$ (Planck's equation).  If their environment has temperature $T$, then the thermal energy per degree of freedom is $kT$ (where $k$ is Boltzmann's constant).  Whether or not the particles behave as quantum particles then depends on two energy scales $h\nu$ and $kT$.  If $kT$ is much less than $h\nu$, then the physics of particles are within the quantum realm.  For a quantum system to be protected from noise, the quantum system needs to be well isolated from the outside world.  Therefore, a particle with characteristic frequency $\nu$ has a \lq\lq quantum temperature" $T_Q=\frac{h\nu}{k}$, which sets an energy scale of the particle.  A modeling as a classical or quantum system is determined by the competition of $T$ and $T_Q$ so that as long as $T_Q$ is much smaller than $T$, then it is safe to regard the system as a quantum system at zero temperature.

Universal properties emerge from the interaction or arrangement of particles at low energy and long wave length.  Since atoms in solids are regular arrays, lattices are real physical objects.  Therefore, we start with a Hilbert space which is a tensor product of a small Hilbert space associated to each atom.  The small Hilbert space comes with a canonical basis and the large Hilbert space has the tensor product basis, which represents classical configurations of the physical system.  Hence, lattices and bases of Hilbert spaces should be taken seriously as physical objects.  Lattices in the following discussion will be triangulations of space manifolds.  In TPMs, we are interested in theories that have continuous limits, i.e. are independent of triangulations.  Such theories are rare.

\subsection{Physical Quantization}

Where do Hilbert spaces come from?  Where do the model Hamiltonians and observables come from?  What are the principles to follow? What is a phase diagram? What is a phase transition? What does it mean that a model is exactly solvable or rigorously solvable?
We cannot answer all these questions, and when we can our subject will be mature and cease to be exciting.

Physicists have a powerful conceptual method to quantize a classical observable $o$: put a hat on the classical quantity, which turns the classical observable $o$ into an operator $\hat{o}$.  For example, to quantize a classical particle with position $x$ and momentum $p$, we define two Hermitian operators $\hat{x}$ and $\hat{p}$, which satisfy the famous commutator $[\hat{x},\hat{p}]=i\hbar$.  The Hilbert space of states are spanned by eigenstates of $\hat{x}$ denoted as $|x\rangle$, where $x$ is the corresponding eigenvalue of the Hermitian operator $\hat{x}$.  Similarly, we can use the momentum $p$, but not both due to the uncertainty principle.  The two quantizations using $x$ or $p$ are related by the Fourier transform.  Therefore, physical quantization is easy: wear a hat.  In general, a wave function $\Psi(x)$ is really a spectral function for the operator $\Psi(\hat{x})$ at state $|\Psi(x)\rangle$: $\Psi(\hat{x})|\Psi(x)\rangle=\Psi(x)|\Psi(x)\rangle$.

\subsection{Phases of Matter}

Relativity unveils the origin of ordinary matter. Condensed matter come in various phases\footnote{The words state and phase are used interchangeably for states of mater. Since we also refer to a wave function as a quantum state, we will use phase more often.}: solid, liquid, and gas.  By a more refined classification, each phase consists of many different ones.  For example, within the crystalline solid phase, there are many different crystals distinguished by their different {\it lattice}\footnote{Lattices in the physical sense: they are regular graphs, not lattices in the mathematical sense that they are necessarily subgroups of $\mathbb{R}^n$ for some $n$.} structures.  Solid, liquid, and gas are classical in the sense they are determined by the temperature. More mysterious and challenging to understand are quantum states of matter: phases of matter at zero temperature (in reality very close to zero).  The modeling and classification of quantum phases of matter is an exciting current research area in condensed matter physics and TQC.  In recent years, much progress has been made in a particular subfield: TPMs.  Besides their intrinsic scientific merits, another motivation comes from the potential realization of fault-tolerant quantum computation using non-abelian anyons.

Roughly speaking, a phase of matter is an equivalence class of quantum systems sharing certain properties.  The subtlety is in the definition of the equivalence relation.  Homotopy class of Hamiltonians is a good example. Phases are organized by the so-called phase diagram: a diagram represents all possible phases and some domain walls indicating the phases transitions.  Quantum systems in the same domain are in the same phase, and the domain walls are where certain physical quantities such as the ground state energy per particle become singular.

Matter is made of atoms and their arrangement patterns determine their properties.  One important characteristic of their patterns is their symmetry.  Liquids have a continuous symmetry, but a solid has only a discrete symmetry broken down from the continuous.  In $3$ spatial dimensions, all crystal symmetries are classed into $230$ space groups.
A general theory for classical phases of matter and their phase transitions was formulated by Landau: In his theory,
phases of matter are characterized by their symmetry groups, and phase transitions are characterized by symmetry breaking.
It follows that group theory is an indispensable tool in condensed matter physics. TPMs do not fit into the Landau paradigm, and currently 
intensive effort is being made in physics to develop a post-Landau paradigm to classify quantum phases of matter.

\subsection{Quantum Qudit Liquids}\label{topql}

Two fascinating macroscopic quantum phenomena are high temperature superconductivity and the fractional quantum Hall effect.  Both classes of quantum matter are related to topological qubit liquids\footnote{The common term for them in physics is \lq\lq spin liquids", as we will call them sometimes too.  But the word \lq\lq spin" implicitly implies $\SU(2)$ symmetry for electron systems, which is not present in general. Therefore, we prefer \lq\lq qubit liquids" or really \lq\lq qudit liquids", and use \lq\lq spin liquids" only for those with an $\SU(2)$ symmetry.}.  While the quantum Hall liquids are real examples of topological phases of matter, the role of topology in high temperature superconductors is still controversial.  But the new state of matter \lq\lq quantum spin liquid" suggested for high temperature superconductors by Anderson is possibly realized by a mineral: herbertsmithite.

Herbertsmithite is a mineral with the chemical formula $\textrm{Cu}_3\textrm{Zn}\textrm{(OH)}_6 \textrm{Cl}_2$.  The unit cell has three layers of copper ions $\textrm{Cu}^{2+}$ which form three perfect kagome lattices\footnote{Kagome is a traditional Japanese woven bamboo pattern.  By the Kagome lattice, we mean the graph consisting of the vertices and edges of the trihexagonal tiling alternating triangles and hexagons.} in three parallel planes.  The copper ion planes are separated by zinc and chloride planes.  Since the $\textrm{Cu}^{2+}$ planes are weakly coupled and $\textrm{Cu}^{2+}$ has a $S=\frac{1}{2}$ magnetic moment, herbertsmithite can be modeled as a perfect $S=\frac{1}{2}$ kagome antiferromagnet with perturbations.

\subsubsection{Fundamental Hamiltonian vs model Hamiltonian}

Theoretically, the {\it fundamental}\footnote{By fundamental here, we mean the Hamiltonian comes from first physical principles such as Coulomb's law.  Of course they are also model Hamiltonians philosophically. By model here, we mean that we describe the system with some effective degrees of freedom, and keep only the most relevant part of the interaction and treat everything else as small perturbations.} Hamiltonian for the herbertsmithite can be written down, and then we just need to find its ground state wave function and derive its physical properties.  Unfortunately, this fundamental approach from first principle physics cannot be implemented in most  realistic systems.  Instead an educated guess is made: a model Hamiltonian is proposed, and physical properties are derived from this model system.  This emergent approach has been extremely successful in our understanding of the fractional quantum Hall effect \cite{Wenbook}.

The model Hamiltonian for herbertsmithite is the Heisenberg $S=\frac{1}{2}$ kagome antiferromagnet with perturbations: $$ H=J\sum_{\langle i,j\rangle} S_i\cdot S_j+ H_{\textrm{pert}},$$ where $J$ is the exchange energy, and $S_i$ are the spin operators, i.e. $2\times 2$ Pauli matrices.  An explicit descrition of $H$ is given in Example \ref{Kagome} below and $H_{\textrm{pert}}$ consists of certain $3$-body terms with operator norm much smaller than $J$. Experiments show that $J=170K\sim 190K$.  The indices $i,j$ refer to the vertices of the lattice (also called sites in physics representing the phyiscal copper ions $\textrm{Cu}^{2+}$), and ${\langle i,j\rangle}$ means that the sum is over all pairs of vertices that are nearest neighbors.

There are many perturbations (small effects) to the Heisenberg spin exchange Hamiltonian $H$.  Different perturbations lead to different potential spin liquid states: numerical simulation shows that next nearest neighbor exchange will stabilize a gapped spin liquid, while other models lead to gapless spin liquids.

\subsection{Many-body Quantum Systems}\label{mathqs}

\subsubsection{Linear algebra problems need quantum computers}

To a first approximation, the subject of quantum many-body systems is linear algebra in the quantum theory language.  In practice, it is a linear algebra problem that even the most powerful classical computer cannot solve.  Memory is one major constraint.  State-of-the-art computational physics techniques can handle Hilbert spaces of dimension $\approx 2^{72}$.  Compared to a real quantum system, $2^{72}$ is a small number: in quantum Hall physics, there are about $10^{11}$ electrons per $\textrm{cm}^2$---thousands of electrons per square micron.  The Hilbert space for the electron spins has dimension $\approx 2^{1000}$.  Therefore, it is almost impossible to solve such a problem exactly.  To gain understanding of such problems, we have to rely on ingenious approximations or extrapolations from small numbers of electrons.

\begin{definition}

\begin{enumerate}

\item A many-body quantum system (MQS) is a triple $(\mcL, b, H)$, where $\mcL$ is a Hilbert space with a distinguished orthonormal basis\footnote{The preferred basis defines locality, which is very important for defining entanglement.  The role of basis is also essential for superposition to be meaningful.} $b=\{e_i\}$,  and $H$ a Hermitian matrix regarded as a Hermitian operator on $\mcL$ using $b$.  The Hermitian operator $H$ is called the Hamiltonian of the quantum system, and its eigenvalues are the energy levels of the system. The distinguished basis elements $e_i$ are the initial classical states or configurations.

\item An MQS on a graph $\Gamma=(V,E)$ with local degrees of freedom (LDF)---a qudit space---$\mbbC^d$ (this is the so-called LDF) is an MQS $(\mcL, b,H)$ with Hilbert space $\mcL=\otimes_{e\in E} \mbbC^d$ where $E$ are the edges (bonds or links) of $\Gamma$, the orthonormal basis $b$ is obtained from the standard basis of $\mbbC^d$, and some local Hamiltonian $H$. A Hamiltonian is local if it is $k$-local for some $k$ as in definition 4.5. Interesting Hamiltonians usually result from Hamiltonian schemas defined below.

In quantum computing jargon, there is a qudit on each edge.
We will use the Dirac notation to represent the standard basis of $\mbbC^d$ by $e_i=|i-1\rangle, i=1,\ldots,d$.  When $d=2$, the basis elements of $\mcL$ are in one-one correspondence with bit-strings or $\mbbZ_2$-chains of $\Gamma$.

\end{enumerate}

\end{definition}

Since our interest is in quantum phases of matter, we are not focusing on a single quantum system.  Rather we are interested in a collection of MQSs and their properties in some limit, which corresponds to the physical scaling or low-energy/long-wave length limit.  In most cases, our graphs are the $1$-skeleton of some triangulation of a manifold.

The most important Hermitian matrices are the Pauli matrices, which are spin=$\frac{1}{2}$ operators.  Pauli matrices wear two hats because they are also unitary.

\begin{example}\label{Kagome}

Let $\Gamma$ be the kagome lattice on the torus, i.e., a kagome lattice in the plane with periodic boundary condition.  The Hilbert space consists of a qubit on each vertex and $\langle i, j\rangle$ represents an edge with vertices $i$ and $j$.  The spin=$\frac{1}{2}$ Heisenberg Hamiltonian is

$$H=J \sum_{\langle i,j\rangle} \sigma^x_i\sigma^x_j+\sigma^y_i\sigma^y_j+\sigma^z_i\sigma^z_j,$$
where the Pauli matrix $\sigma_k^\alpha$, for $\alpha=x,y,z$ and $k=i,j$, acts by $\sigma^\alpha$ on the $k$th qubit and by $\textrm{Id}$ on the other factors. The coupling constant $J$ is some positive real number.

\end{example}

\begin{definition}

Let $(\mcL, b, H)$ be an MQS and $\{\lambda_i, i=0,1,\ldots\}$ the eigenvalues of $H$ ordered in an increasing fashion and $\mcL_{\lambda_i}$ the corresponding eigenspace.

\begin{enumerate}

\item $\mcL_{\lambda_0}$ is called the ground state manifold and any state in $\mcL_{\lambda_0}$ is called a ground state.  Any state in its complement $\oplus_{i>0}\mcL_{\lambda_i}$ is called an excited state.  Usually we are only interested in the first excited states in $\mcL_{\lambda_1}$, but sometimes also states in $\mcL_{\lambda_2}$.  Bases states in $\mcL_{\lambda_i}$ are called minimal excited states.

\item  A quantum system is rigorously solvable if the ground state manifold and minimal excitations for low excited states are found.  A quantum system is exactly solvable if the exact answers are known physically (but not necessarily rigorously).

\end{enumerate}

\end{definition}

Anyons are elementary excitations of TPMs, but it is difficult to define the notion mathematically.  Morally, an elementary excitation is minimal if it is not a linear combination of other non-trivial excitations, but due to constraints from symmetries such as in the toric code, the minimal excitations are anyon pairs created out of the vacuum.

\begin{definition}

A unitary $U$ is a symmetry of a Hamiltonian $H$ if $UHU^{\dagger}=H$.  A Hermitian operator $K$ is a symmetry of $H$ if $[K,H]=0$, and then $e^{itK}$ is a unitary symmetry of $H$ for any $t$.  When $U$ is a symmetry, each energy eigenspace $\mcL_{\lambda_i}$ is further decomposed into eigenspaces of $U$.  Those eigenvalues of $U$ are called good quantum numbers.

\end{definition}

\begin{definition}

\begin{enumerate}

\item Expectation values: the expectation value of an observable $\mcO$ is $\langle 0|\mcO|0 \rangle$. Recall an observable in a quantum theory is just a Hermitian operator.

\item Correlation functions:  the correlation functions of $n$ observables $\mcO_i$ at sites $r_i$ are $$\langle 0|\mcO_n(r_n)\cdots \mcO_{1}(r_1)|0 \rangle.$$

\end{enumerate}

\end{definition}

\begin{definition}

\begin{enumerate}

\item Given an MQS on a graph with Hilbert space $\mcL=\bigotimes_{i\in I} L_i$, where the $L_i$ are local Hilbert spaces such as the same qudit, an operator $O$ is $k$-local for some integer $k>0$ if $O$ is of the form $\textrm{Id}\otimes A\otimes \textrm{Id}$, where $A$ acts on a disk of radius=$k$ at a vertex in the edge length metric.

\item An MQS is $k$-local for some integer $k>0$ if $H=\sum_i H_i$ with each $H_i$ $k$-local.

\item A Hamiltonian is a sum of commuting local projectors (CLP) if $H=\sum_i H_i$ such that each $H_i$ is a local projector and $[H_i,H_j]=0$ for all $i,j$.

\end{enumerate}

\end{definition}

\subsection{Hamiltonian Definition of Topological Phases of Matter}

TPMs are phases of matter whose low energy universal physics is modeled by
TQFTs and their enrichments. To characterize TPMs, we focus
either on the ground states on all space surfaces or their elementary excited states in the plane. The ground state dependence on the topology of spaces organizes into UTMFs, and elementary excitations as UMCs.

\subsubsection{Pachner poset of triangulations}

\begin{figure}[h]\includegraphics[width=2.5in]{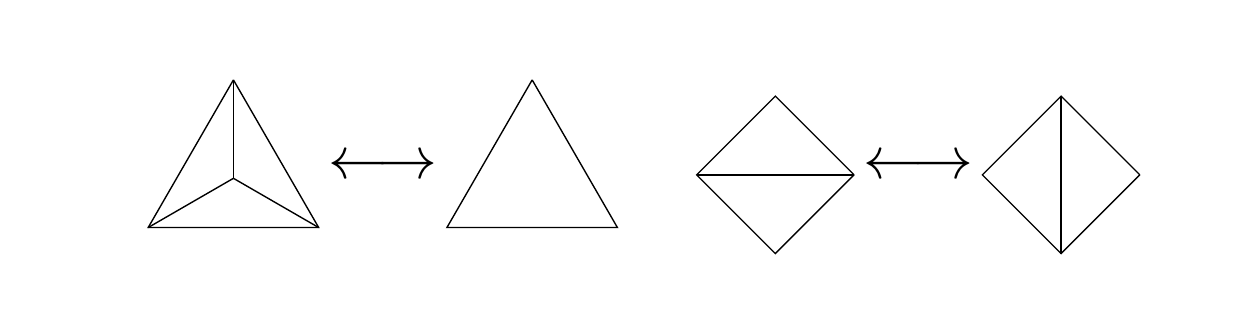}
\caption{\label{Pachner}Pachner moves.}
\end{figure}

Let $Y$ be a closed oriented surface, and $\Delta_i,i=1,2$  two triangulations of $Y$.  By the Pachner theorem, $\Delta_1$ to $\Delta_2$ can be transformed to each other via two type of moves, see Figure \ref{Pachner}.  For triangulations $\Delta_i, i=1,2$ we say $\Delta_1 \leq \Delta_2$ if $\Delta_1$ is a refinement of $\Delta_2$ via finitely many Pachner moves.  Then the partial order $\leq $ makes all triangulations of $\Delta_i$ of $Y$ into a poset.  

\subsubsection{Hamiltonian schemas}

To define a Hamiltonian is to give rules for specifying a class of quantum systems.  Such rules for defining quantum systems will be called a Hamiltonian schema.  In this section, Hamiltonian schemas are given for triangulations of, or lattices in, closed surfaces.  Usually, we also need additional structures on the surfaces or lattices.  In real materials, the lattice points are atoms that form special lattices.  Here we regard every triangulation of a surface as a lattice, therefore our theories are very special: they exist on every space surface with an arbitrary triangulation.  This is a strong locality requirement as such a theory is fully determined by the Hamiltonian on a disk region.

\begin{definition}

A Hamiltonian schema is a uniform local rule to associate an MQS to any triangulation of a surface, where local means the rule is determined by the regular neighborhood of one simplex of each dimension, and uniform means all resulting Hamiltonians are $k$-local for some fixed integer $k>0$.

\end{definition}

\subsubsection{Scaling limits}

Given a closed surface $Y$ and a rule $F$ to assign a vector space $V_F(Y,\Delta_i)$ to each triangulation $\Delta_i$ of $Y$.  Then $V_F(Y,\Delta_i)$ form an inverse system over the index poset of triangulations of $Y$.  The scaling limit $V_F(Y)$ of the rule $F$ is the inverse limit of $V_F(Y,\Delta_i)$ if it exists.

\subsubsection{Energy gap}

\begin{definition}

A Hamiltonian schema $H$ is gapped or has an energy gap if the spectrum of every resulting Hamiltonian has the following structure:  for each closed surface $Y$, there exist some constants $\Lambda >0$ and $\xi>0$, and a set of eigenvalues of $\{\lambda_0^{(j,a)},j=1,2,\ldots,k\}$ of the resulting Hamiltonian $H(Y;\Delta_a)$ from each triangulation $\Delta_a$ such that $|\lambda_0^{(i,a)}-\lambda_0^{(j,a)}|\leq e^{-\xi c(\Delta_a)}$ and all other eigenvalues $\{\lambda^{(a)}_l\}$ satisfying $|\lambda^{(a)}_l-\lambda_0^{(j,a)}| \geq \Lambda$, for all $l\geq 1,i,j=1,2,\ldots,k$.

\end{definition}

In the thermodynamic limit, the energy scales $\{\lambda_0^{(j,a)},j=1,2,\ldots,k\}$ converge to the same low energy level, which is separated from higher energies by the positive constant $\Lambda$.  The existence of such an energy gap has deep implications for the physical systems.  
Establishing a gap of Hamiltonian schemas is a difficult problem.  One obvious case is for Hamiltonians which are CLPs.

\begin{prop}

If $H=\sum H_i$ is a CLP then it has an energy gap. If each $H_i$ is positive semi-definite, then a state $\Psi$ is a ground state if and only if $H_i \Psi=0$ for each $i$.

\end{prop}

\subsubsection{Ground state modular functor}

A gapped Hamiltonian schema $H$ for closed surfaces has a scaling limit if the eigenspaces $\bigoplus_{j=1}^k V_{\lambda_0^{(j,a)}}$ have a scaling limit $V_H(Y)$ for each closed surface $Y$.

\begin{definition}

A gapped Hamiltonian schema $H$ for closed surfaces $Y$ with scaling limit is topological if the ground state functor $Y\rightarrow V_H(Y)$ is naturally isomorphic to a UTMF $V$ on closed surfaces.  Two topological Hamiltonian schemas are equivalent if their ground state functors are equivalent as tensor functors.

\end{definition}

\subsubsection{Stability and Connectivity}

\begin{definition}

\begin{enumerate}

\item A topological Hamiltonian schema is stable if all resulting Hamiltonians satisfy the error-correction TQO1 and ground state homogeneity TQO2 axioms in \cite{BHM10}.

\item Two stable gapped Hamiltonian schemas $H$ and $H'$ are connected if there exists a path of gapped Hamiltonian schemas $H_t, 0\leq t\leq 1$ such that $H_0=H$ and $H_1=H'$.

\end{enumerate}

\end{definition}

A stable Hamiltonian represents a phase of matter as small perturbations will not change the physical properties of the quantum system, so they all represent the same phase.  TQO1 and TQO2 are  conditions sufficient to guarantee stability of gapped Hamiltonians.  It is not impossible that they are automatically satisfied for all topological Hamiltonian schemas.

\begin{definition}

A $2D$ TPM is an equivalence class of connected stable topological Hamiltonian schemas.

The {\it topological order} in a TPM is the UMC derived from the ground state UTMF.  This undefined physical notion conveys a vague picture of a many particle system with a dynamical pattern.  Here we use it to characterize topological phases of matter so that two topological phases of matter represent the same phase if their topological orders are the same.

\end{definition}

\begin{conjecture}

Elementary excitations of a TPM form a UMC, and connected stable Hamiltonian schemas are equivalent.

\end{conjecture}

\subsection{Realization of UTMF}

Given a UTMF $V$, can $V$ be realized by a TPM?

\begin{conjecture}
A UTMF $V$ can be realized by a TPM if and only if its associated UMC $\mcB$ is a quantum double/Drinfeld center $\mcZ(\mcC)$ of a unitary fusion category $\mcC$.
\end{conjecture}

Indeed there is no doubt that any UTMF with a doubled UMC can be realized by the Levin-Wen model (LW) \cite{LW05}, but a fully rigorous mathematical proof in our sense is still not in the literature.  A step towards the conjecture is made in \cite{kirillov11}.  To complete a proof in our sense, the string-net space needs to be identified with the exact ground states of the microscopic Hamiltonians of the LW model in a functorial way, and the Drinfeld center $\mcZ(\mcC)$ with the excited states.

\subsubsection{Hamiltonian realization of Dijkgraaf-Witten theories}

The untwisted Dijkgraaf-Witten TQFTs are constructed from finite groups $G$ \cite{DijWit}.  Kitaev introduced Hamiltonian schemas to realize them \cite{Ki03}, and the $G=\mathbb{Z}_2$ case is the celebrated toric code.

\subsubsection{Hamiltonian realization of Turaev-Viro theories}

The LW Hamiltonian schemas, which generalizes Kitaev schemas dually, should realize Barrett-Westbury-Turaev-Viro TQFTs.  
As input, the LW model takes a unitary fusion category $\mathcal{C}$.
For example, there are two natural unitary fusion categories associated to a finite group $G$:
\begin{enumerate}
\item The category $\textrm{Vec}_G$ of $G$=graded vector spaces, with simple objects labeled by the elements of $G$.
\item The representation category $\textrm{Rep}(G)$ of $G$ with simple objects the irreducible representations of $G$.
\end{enumerate}
These two categories are monoidally inequivalent, but Morita equivalent.  The LW model with these inputs realize the same Drinfeld center $\mcZ(\textrm{Vec}_G)=\mcZ(\textrm{Rep}(G))$. 

\subsubsection{Ground states and error correction codes}

TPMs are natural error-correcting codes.
Indeed the disk axiom of UTMFs implies that local errors are scalars: if an operator is supported on a disk, then splitting the disk off induces a decomposition of  the modular functor space $V(\Sigma) \cong V(\Sigma') \otimes V(D^2)$, where $V(D^2) \cong \mbbC$ and $\Sigma'$ is the punctured surface. So TQO1 of \cite{BHM10} is basically the disk axiom.  The following is true for the toric code \cite{Ki03}:

\begin{conjecture}
The ground states of the LW model form error-correcting codes.
\end{conjecture}

\subsection{Realization of UMC}

Can a UMC $\mcB$ always be realized by elementary excitations of a TPM?
In principle, the UMC $\mcB$ should be derived from elementary excitations in the plane with string operators.  There is a large physics literature, but a rigorous mathematical proof is still highly non-trivial.

\subsection{Other Definitions of Topological Phases of Matter}

Besides the Hamiltonian definition of TPMs, another common definition in physics 
of TPM is via states.  A state is “trivial” if it can be changed to 
a product state using a local quantum circuit. Otherwise the state is topologically ordered \cite{Chen12}. When defined properly, all definitions should be equivalent to each other, but mathematically the problem is fundamentally unclear.

\subsection{Open Problems}

A UMC $\mcB$ is chiral if the central charge of $\mcB$ is not $0$ mod $8$, i.e. $\sum_{i\in \Pi_\mcB} d_i^2 \theta_i \neq \sqrt{\dim(\mcB)}$.  A chiral TPM is one with a chiral topological order, i.e. a chiral UMC.

One of the most difficult problems is to realize chiral topological orders by Hamiltonian schemas.  Numerical simulation undoubtedly suggests the two examples below realize chiral TPMs.  But since the Hamiltonians are not CLPs, mathematical proofs seem to be completely out of reach. 

\begin{conjecture}

There are no CLP Hamiltonian realizations of chiral TPMs.

\end{conjecture}

\subsubsection{Haldane Hamiltonian for semions}

Numerical simulation \cite{vidalsemion} shows that Haldane Hamiltonian in \cite{haldane88} realizes the semion theory for certain coupling constants.
Proving this mathematically would provide the first mathematical example of a chiral TPM in our sense.

\subsubsection{Kitaev model for Ising}

The Kitaev honeycomb model \cite{Ki06} should realize mathematically the Ising modular category.  

\section{Anyonic Quantum Computation}\label{Sec:AQC}

\lq\lq {\textit{Quantum computation} is any computational model based upon
the theoretical ability to manufacture, manipulate and measure
quantum states''\cite{FKLW03}.  In TQC, information is encoded in multi-anyon quantum states, and our goal is the construction of a large scale quantum computer based on braiding non-abelian anyons.  In this section, we cover the anyonic quantum computing models using anyon language.  

To carry out quantum computation, we need quantum memories, quantum gates, and protocols to write and read information to and from the quantum systems.  In the anyonic quantum computing model,  we first pick a non-abelian anyon type, say $x$.  Then information is stored in the ground state manifold $V_{n,x;t}$ of $n$ type $x$ anyons in the disk with total charge $t$.  As $n$ goes to infinity, the dimension of $V_{n,x;t}$ goes asymptotically as $d_x^n$, where $d_x>1$ is the quantum dimension of the non-abelian anyon $x$.  It follows that when $n$ is large enough, we can encode any number of qubits into some $V_{n,x;t}$.  The ground state manifold $V_{n,x;t}$ is also a non-trivial unitary representation of the $n$-strand braid group $\mcB_n$ \cite{RWpra}, and the unitary representation matrices serve as quantum circuits.  An initial state in a computation is given by creating anyons from the ground state and measurement is achieved by fusing anyons together to observe the possible outcomes.  Subtleties arise for encoding qubits into $V_{n,x;t}$ because their dimensions are rarely powers of fixed integers.  Another important question is whether the braiding matrices alone will give rise to a universal gate set.

The weak coupling of topological degrees of freedom with local ones is both a blessing and curse.  The topological protection is derived directly from this decoupling of topological degrees of freedom from the environment.  On the other hand, the decoupling also makes  experimental detection of topological invariants difficult as most experiments measure local quantities such as electric currents.

\subsection{Topological Qudits}

There are many choices to encode strings $b\in \mbbZ_d^\ell$ onto topological degrees of freedom in ground states $V_{n,x;t}$ of $n$ anyons of type $x$.  The explicit topological encoding uses the so-called fusion channels of many anyons, so that strings $b$ correspond to fusion-tree basis elements \cite{Wa10}.\footnote{Earlier encoding in \cite{OP99} involves a splitting of certain fusion channels and a reference bureau of standards, hence not explicitly topological.  Though the computation might still be carried out using only the topological degrees of freedom, the relation of the two encodings is not completely understood, and not analyzed carefully in the literature.  Our encoding follows the explicitly topological one in \cite{FKLW03}.}

\subsubsection{Dense vs sparse encoding}

The sparse encoding is directly modeled on the quantum circuit model, so topological subspaces are separated into single qudits.  In the dense encoding, qudits are encoded into topological subspaces, but no separation into qudits is provided.  There are protocols to go from one to the other by using measurements, and they are not equivalent in general.  We will only use sparse encoding in the following.

The sparse encoding of one and two qudits by fusion trees is shown in Figure \ref{sparseencoding}, where the labels $a_i,b_j, t_k$ are anyons resulting from the fusions.
The encoding for $m$-qudits is analogous.

\begin{figure}[hbp]
\centering
\begin{picture}(160,60)(-20,-10)
 \label{1-qudit pic}
\put(50,0){\line(0,1){10}}
\put(50,10){\line(1,1){30}}
\put(50,10){\line(-1,1){30}}
\put(70,30){\line(-1,1){10}}
\put(30,30){\line(1,1){10}}

\put(20,42){$X$}
\put(40,42){$X$}
\put(60,42){$X$}
\put(80,42){$X$}
\put(30,19){$a$}
\put(66,19){$b$}
\put(51,-2){$t_1$}
\end{picture}
\begin{picture}(180,120)(0,-50)
\put(50,10){\line(1,-1){45}}
\put(50,10){\line(1,1){30}}
\put(50,10){\line(-1,1){30}}
\put(70,30){\line(-1,1){10}}
\put(30,30){\line(1,1){10}}

\put(140,10){\line(-1,-1){45}}
\put(140,10){\line(1,1){30}}
\put(140,10){\line(-1,1){30}}
\put(160,30){\line(-1,1){10}}
\put(120,30){\line(1,1){10}}

\put(95,-35){\line(0,-1){15}}

\put(20,42){$X$}
\put(40,42){$X$}
\put(60,42){$X$}
\put(80,42){$X$}
\put(30,18){$a_1$}
\put(65,18){$b_1$}
\put(58,-15){$t_1$}

\put(110,42){$X$}
\put(130,42){$X$}
\put(150,42){$X$}
\put(170,42){$X$}
\put(120,18){$a_2$}
\put(154,18){$b_2$}
\put(126,-15){$t_1$}

\put(86,-47){$t_2$}
\end{picture}
\caption{\label{sparseencoding}(Left) one qudit and (Right) two qudits.}
\end{figure}
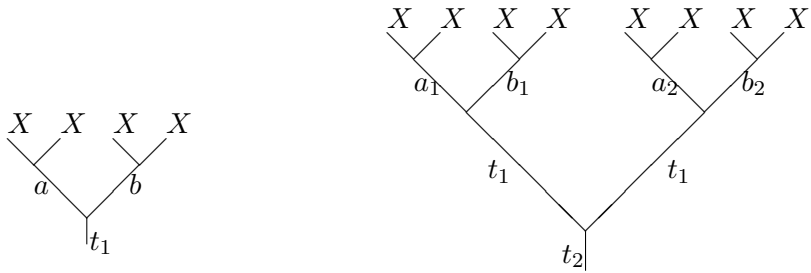

\subsection{Topological Gates}

Straightforward topological gates are braiding gates.  But there are also topological resources that can be used to augment braiding gates such as measurements and mapping class group representations of higher genus surfaces.  We will mainly focus on braiding gates and the simplest measurement:  the measurement of total charge of a group of anyons.

\subsubsection{Braiding gates and universality}

Given an anyon $X$ and $t_i, i=1,2$, $V_{4,X;t_1}=\textrm{Hom}(X^{\otimes 4},t_1)$ and $V_{8,X;t_2}=\textrm{Hom}(X^{\otimes 8},t_2)$ are representations of the braid groups $\mcB_4$ and $\mcB_8$, respectively.  Labeled trees as in Figure \ref{sparseencoding} form orthonormal bases of these representations, and $V_{8,X;t_1,t_1,t_2}$ is a subspace of $V_{8,X;t_2}$ consisting of a subset of the labeled trees, which is not necessarily a sub-representation of $\mcB_8$.  With respect to these tree bases, $\rho_X (\sigma_i)$ are unitary matrices, which are quantum gates for quantum computing.

\begin{definition}

Given an anyon $X$ and an encoding of one qudit and two qudits as above, the braiding gate set from the anyon $X$ consists of the $1$-qudit gates $\rho_X (\sigma_i), i=1,2,3$, and $2$-qudit gates $\rho_X (\sigma_i), i=1,2,...,7$ that preserve the qudits $V_{4,X;t_1}$ and $V_{8,X;t_1,t_1,t_2}$, respectively.

\end{definition}

In practice, the above definition of braiding gates is too restrictive as gates $\rho_X (\sigma_i)$ such as above without leakage is rare.  In general, the braid gate gate $\rho_X (\sigma_i)$ usually has leakage in the sense that it does not preserve the qudits $V_{4,X;t_1}$ and $V_{8,X;t_1,t_1,t_2}$.  The general phenomenon is that leakage-free gates are weak for quantum computation in the sense we do not know any leakage free universal braiding gate set.  This is closely related to the Property $F$ Conjecture \ref{propF} and Conjecture \ref{localizeconj}.

To achieve universality, we have two options: for leakage-free gate sets, we supplement them with measurements, and for those without leakage-free options, we prove density, and then leakage can be basically eliminated.

One might wonder if it is possible to simulate the braiding gates for a TQC via quantum circuits without a large non-computational component.  This is a kind of converse to ``leakage" described in \cite{Wa10}, which leads to the notion of \emph{localization} \cite{RWlocal,GHW13}.  

\begin{definition}\label{localize}
A \emph{localization} of an anyon $X$ is a unitary braided vector space $(R,W)$ (i.e. a unitary matrix $R\in\Aut(W^{\otimes 2})$ that satisfies the Yang-Baxter equation) and injective algebra maps $\tau_n$ so that the following diagram commutes for all $n$: $$\xymatrix{ \mbbC\mcB_n\ar[d]^{\rho_X}\ar[dr]^{\rho^R} \\
\mbbC\rho_X(\mcB_{n})\ar[r]^{\tau_n} & \End(W^{\otimes n})}$$ where $\rho_X:\mcB_n\rightarrow \Aut(\bigoplus_{a\in\Pi_\mcC}\Hom(a,X^{\otimes n}))$ is the corresponding braid group representation, $\rho^R$ the braid group representation from the $R$-matix $R$, and $\mbbC\rho_X(\mcB_{n})$ is the group algebra of the image $\rho_X(\mcB_{n})$. When such a localization exists the braiding gates $\rho_X(\sigma_i)$ can be simulated on quantum circuits using the gate $R$.
\end{definition}

\begin{definition}\label{BraidingUniversal}

An anyon $X$ is called \emph{braiding universal} if, for some $n_0$, the images of $\mcB_n$ on the irreducible sub-representations $V\subset\End(X^{\otimes n})$ are dense in $SU(V)$ for all $n\geq n_0$.

\end{definition}

There is a small subtlety related to whether the braid group images are independently dense on each $V$ \cite{Wa10}, but we will not discuss the  details for making a computing model out of a braiding universal anyon here.  We illustrate the general idea as below and refer the interested readers to \cite{FLW02,Wa10}.

To perform a quantum circuit $U_L$, we want to find a braid $b$ such that the braiding matrix $U(b)$ from the representation of $b$ makes the following diagram commute:

$$\xymatrix{
{{(\mathbb{C}^2)}^{\otimes n}} \ar[d]^{U_L} \ar[r] & V(Y) \ar[d]^{U(b)}\\
{{(\mathbb{C}^2)}^{\otimes n}} \ar[r] & V(Y) }$$

For most anyons, this diagram does not commute exactly.  So we need to work with approximations, which is sufficient theoretically.

\subsubsection{Resource-assisted universality}

Besides braidings, there are other topological operations such as measuring total topological charges \cite{Bonderson08}, using ancillary topological states \cite{Bravyi05}, using symmetry defects to access mapping class group representations \cite{BF16}, and gapped boundaries \cite{CMW17,CMW16}.

Braiding gate sets can be supplemented by gates from those topological operations, thus to obtain measurement assisted, ancillary assisted, or genus-assisted universal gate set \cite{Ki03,Bonderson08,Cui151,Cui152,Bravyi00,Bonderson10,CMW16}.

\subsection{Density of TQFT Representations}

Definition \ref{BraidingUniversal} raises a purely mathematical question: How can we detect whether or not a given anyon is braiding universal?  The Property $F$ Conjecture \ref{propF} asserts that an anyon $X$ is  braiding universal if and only $d_x^2\not\in \mathbb{Z}$ (\cite{NR11}).

\subsubsection{The N-eigenvalue problem}

Given a particular anyon type $x$, we analyze the braid group representation as follows:
\begin{enumerate}
\item Determine if the braid representations $V_{n,x;t}=\Hom(t,x^{\otimes n})$ are irreducible for all $n$.
This turns out to be a very difficult question in general.
If reducible, we must decompose it into irreducible representations (irreps).
\item The number of distinct eigenvalues of the braiding $c_{x,x}$ is bounded by $\sum_{i \in L} N^i_{x,x}$.
Since all braid generators are mutually conjugate, the closure of $\rho_{n,x;t}(\mcB_n)$ in $\mbbU(V_{n,x;t})$ is generated by a single conjugacy class.
\end{enumerate}

\begin{definition}
Let $N \in \mathbb{Z}_+$. We say a pair $(G,V)$, $G$ a compact Lie group, $V$ a faithful irrep of $G$, has the $N$-eigenvalue property if there exists an element $g \in G$ such that the conjugacy class of $g$ generates $G$ topologically and the spectrum $X$ of $\rho(g)$ has $N$ elements and satisfies the no-cycle property: $u \{ 1, \xi, \xi^2, \ldots, \xi^{n-1} \} \not \subset X$ for any $n$th root of unity $\xi$, $n \geq 2$, and $u \in \mathbb{C}^\times$.
\end{definition}

The $N$-eigenvalue problem is to classify all pairs with the $N$-eigenvalue property.
For $N=2,3$, this is completed in \cite{FLW02} and \cite{LRW05}.
As a direct corollary, we have
\begin{theorem} Suppose $r \neq 1,2,3,4,6$, and $n \geq 4$ if $r = 10$.
\begin{enumerate}
\item The closed images of the Jones representations of $\mcB_n$ on $V^a_{n,\frac{1}{2};t}$ contain $\SU(V^a_{n,\frac{1}{2};t})$.
\item The anyons of type $s=\frac{1}{2}$ in $\SU(2)_k$ are braiding universal.
\end{enumerate}
\end{theorem}

It follows that Theorem \ref{jonesbqp} can be strengthened to say ``BQP-complete" for $r\neq 3,4,6$.

\subsection{Topological Quantum Compiling}

\subsubsection{Exact braiding gates}

Given a braid representation $V$ with a fusion tree basis \cite{Wa10}, the matrices in $\mathbb{U}(N)$ in the image of the representation are topological circuits if we use the images of elementary braids as braiding gates.  The fusion tree basis is physical so the resulting unitary matrices are reasonable gates.  It is difficult to decide which circuits can be realized by braids in a given representation.  The only complete answer is for the one-qubit gates for the Fibonacci anyon \cite{KBS14}. 

\subsubsection{Approximation braiding gates}

Exact entangling circuits are extremely difficult to find, so the next question is how to efficiently approximate entangling circuits by braiding ones.

\subsection{Open Problems}

\subsubsection{Exact entangled gates}

In the sparse encoding, we do not known any exact two qudits braiding entangling gates.  In the dense encoding, CNOT can be realized for Ising anyons, but we do not know if there is a leakage-free entangling gate for Fibonacci anyon.

\subsubsection{Computing model from CFTs}

It is interesting to understand what are the natural computing models from CFTs \cite{Moji18}.

\section{On Modular (Tensor) Categories}\label{section: Modular Categories}

All linear categories and vector spaces in this section are over the complex numbers $\mbbC$.  The axiomatic definition of modular categories was given in subsection \ref{sssection: Modular}.

Modular tensor categories first appeared as a collection of tensors in the study of CFTs \cite{MS89}, while the definition of a modular category was formulated to algebraically encode Reshetikhin-Turaev TQFTs \cite{Turaev92}.  Modular tensor categories and modular categories define equivalent algebraic structures \cite{DHW13}, and we use the latter term.  Modular categories arise naturally in a variety of mathematical subjects, typically as representation categories of algebraic structures such as quantum groups/Hopf algebras \cite{BKi}, vertex operator algebras \cite{Huang05}, local conformal nets \cite{KLM}, loop groups and von Neumann algebras \cite{DK98}

Significant progress on the classification of modular categories has been made during the last decade \cite{M2,RSW09,BNRW151,BNRW152}, and a structure theory for modular categories is within reach.  A fruitful analogy is to regard modular categories as quantized finite abelian groups, and more generally spherical fusion categories as quantized finite groups.  A central theme in modular category theory is the extension of classical results in group theory to modular categories such as the Cauchy \cite{BNRW151} and Landau \cite{ENO1} theorems.

A complete classification of modular categories includes a classification of finite groups in the following sense: any finite group $G$ can be reconstructed from its (symmetric fusion) representation category $\textrm{Rep}(G)$ \cite[Theorem 3.2]{DeligneMilne}.  The Drinfeld center $\mcZ(\textrm{Rep}(G))$ of $\textrm{Rep}(G)$ is always a modular category with $\textrm{Rep}(G)$ as a (symmetric) subcategory.  Therefore, the complete classification of modular categories is extremely difficult without certain restrictions.  One possibility is to classify modular categories modulo the classification of Drinfeld centers of finite groups.

Our interest in modular categories come from their applications to TQC and TPMs, where unitary modular categories model anyon systems.  Modular categories form part of the mathematical foundations of TQC and TPMs \cite{NSSFD08,Wa10}, and their classification would provide a sort of ``periodic table'' of these phases of matter.  More generally, many important practical and theoretical questions in TQC and TPMs can be translated into mathematical questions and conjectures for modular categories. The book \cite{etingof15} is an excellent reference for the background materials, and the survey \cite{Mueger12} covered many of the earlier results.

Examples of modular categories appear naturally in the study of representations of finite groups, Hopf algebras/quantum groups and skein theory for quantum invariants.  These ubiquitous examples are typically related to well-known TQFTs and the computational complexity of their corresponding quantum link invariants provided motivation for the subject of TQC.  Just as the heavier elements of the periodic table are synthesized in the lab from naturally occurring elements, new modular categories are constructed from known categories using various tools such as de-equivariantization, Drinfeld centers and gauging.  The general landscape of modular categories is still largely unexplored, but a few recent breakthroughs gives us hope.

\subsection{Basic Examples of Modular Categories}

Well-known modular categories appear naturally as representations of some sufficiently well-behaved algebraic structure.  Most of these can be derived from finite groups or quantum groups.

\subsubsection{Pointed modular categories}

The simplest examples of modular categories are constructed from finite abelian groups with non-degenerate quadratic forms.  A pointed modular category is one that every simple object is invertible.  In fact, every pointed modular category is constructed in this way.

Let $G$ be a finite abelian group.  A function $q: G\rightarrow \mathbb{U}(1)$ is a quadratic form if 1) $q(-g)=q(g)$, and 2) the symmetric function $s(g,h)=\frac{q(g+h)}{q(g)q(h)}$ is bi-multiplicative.  The set of quadratic forms $Q(G)$ on $G$ form a group under point-wise multiplication.  A quadratic form $q$ is non-degenerate if its induced bilinear form $s(g,h)$ is non-degenerate.

The label set $\Pi_\mcC$ of a pointed modular category $\mcC$ is a finite abelian group under tensor product.  A pointed modular category $\mcC$ is fully determined by the pair $(\Pi_\mcC, \theta)$, where $\theta (a)$ is the topological twist of the label $a$.  Therefore, pointed modular categories are simply finite abelian groups endowed with non-degenerate quadratic forms.  The associativity isomorphism $\omega$ and the braiding $b$ of $\mcC$ from $(\Pi_\mcC, \theta)$ is provided by the cohomology class $(\omega, b)$ in the Eilenberg-Maclane third abelian cohomology group $H^3_{ab}(G,\mbbU(1))$, which is isomorphic to the group $Q(G)$ of non-degenerate quadratic forms on $G$ \cite{EMcL}.

\subsubsection{Quantum doubles of finite groups}
From an arbitrary finite group $G$ and a $3$-cocycle $\omega:G\times G\times G\rightarrow \mbbC^{\times}$ one constructs the (twisted) \emph{quantum double} $D^\omega G$ \cite{Dr3}, a quasi-triangular quasi-Hopf algebra with underlying vector space $D^\omega G=(\mbbC G)^*\otimes \mbbC G$.  An integral modular category is one that every simple object has an integral quantum dimension.  The category $\Rep(D^\omega G)$ is an integral modular category and the corresponding TQFTs are the Dijkgraaf-Witten theories \cite{DijWit}, while the associated link invariant essentially counts homomorphisms from the fundamental group of the link complement to the group $G$ \cite{FreedQuinn}.  A modular category $\mcC$ is called \emph{group-theoretical} if $\mcC\subset \Rep(D^\omega G)$.

\subsubsection{Quantum groups, conformal field theories and skein theories}

From any simple Lie algebra $\mathfrak{g}$ and $q\in\mbbC$ with $q^2$ a
primitive $\ell$th root of unity one can construct a ribbon fusion
category $\mcC(\mathfrak{g},q,\ell)$ (see \cite{BKi}).  One may
similarly use semisimple $\mathfrak{g}$, but the
resulting category is easily seen to be a direct product of those
constructed from simple $\mathfrak{g}$.  We shall say these categories (or
their direct products) are of \emph{quantum group type}. There
is an oft-overlooked subtlety concerning the degree $\ell$ of
$q^2$ and the unitarizability of $\mcC(\mathfrak{g},q,\ell)$.
Let $m$ be the maximal number of edges between any two nodes
of the Dynkin diagram for $\mathfrak{g}$ with $\mathfrak{g}$ simple, so that $m=1$ for
Lie types $A$,$D$,$E$; $m=2$ for Lie types $B$,$C$,$F_4$; and $m=3$ for Lie
type $G_2$.
\begin{theorem}
If $m \mid \ell$, then $\mcC(\mathfrak{g},q,\ell)$ is a unitary modular category for $q=e^{\pm \pi i/\ell}$.
\end{theorem}
This theorem culminates a long string of works in the theory of quantum groups, see \cite{Rowell06} for references.
If $m \nmid \ell$, there is usually no choice of $q$ making $\mcC(\mathfrak{g},q,\ell)$ unitary.

In \cite{finkel} it is shown that the tensor category associated with level $k$
representations of the affine Kac-Moody algebra $\hat{\mathfrak{g}}$ is
tensor equivalent to $\mcC(\mathfrak{g},q,\ell)$ for $\ell=m(k+\check{h}_\mathfrak{g})$ where
$\check{h}_\mathfrak{g}$ is the dual Coxeter number.
The central charge of the corresponding Wess-Zumino-Witten CFT is $\tfrac{k\dim \mathfrak{g}}{k+\check{h}\mathfrak{g}}$ \cite{yellowbook}.  For integer levels (i.e. when $m\mid\ell$)
we will use the
abbreviated notation $G_k$ to denote the modular category $\mcC(\mathfrak{g},q,\ell)$ with $q=e^{\pi i/\ell}$ and $\ell=m(k+\check{h}_\mathfrak{g})$. The $\mcC(\mathfrak{g},q,\ell)$ TQFTs are the Reshetikhin-Turaev theories mathematically, and the Witten-Chern-Simons TQFTs physically \cite{Tu94}.  Well-known link invariants are associated with these TQFTs, such as the Jones polynomial (for $SU(2)_k$), the HOMFLYPT polynomial (for $SU(N)_k$) and the Kauffman polynomial (for $SO(N)_k$ and $Sp(N)_k$).

Two well-known mathematical axiomatizations of chiral CFTs ($\chi$CFT) are vertex operator algebras (VOAs) and local conformal nets, which are conjecturally equivalent to each other \cite{CK15,Tener16}. Two fundamental theorems prove that the representation categories of VOAs or local conformal nets with certain conditions are modular categories \cite{KLM,Huang05}.  The best-known examples of modular categories coming from $\chi$CFTs are essentially the same as those coming from quantum groups, since they take affine Kac-Moody algebras as input.

Versions of quantum group type modular categories $SU(2)_k$ can be constructed using Temperley-Lieb-Jones skein theories \cite[Chapter XII]{Tu94} and general skein theories for other modular category $G_k$ \cite{TurWenzl}.  Briefly, the idea is as follows: 1) start with a link invariant $\mcG$ (e.g. Jones, HOMFLYPT or Kauffman polynomials) that admits a functorial extension to (the ribbon category of) tangles $\mathscr{T}$ 2)  use $\mcG$ to produce a trace $\mathrm{tr}_\mcG$ on $\mathscr{T}$ and 3) take the quotient $\overline{\mathscr{T}}$ of $\mathscr{T}$ by the tensor ideal of negligible morphisms, which is essentially the radical of the trace $\mathrm{tr}_\mcG$.  The resulting categories are always ribbon categories (sometimes modular) and may differ from the quantum group type categories is subtle ways, such as in Frobenius-Schur indicators (see below).

\subsection{New Modular Categories from Old}
The easiest way in which a new modular category can be constructed from two given modular categories $\mcC$ and $\mcD$ is via the Deligne (direct) product: $\mcC\boxtimes\mcD$, whose objects and morphisms are just ordered pairs, extended bilinearly.
\subsubsection{Drinfeld center}

A categorical generalization of the quantum double construction described above applies to strict monoidal categories, known as the Drinfeld center.  In the following definition, we will write $x\otimes y$ as $xy$ for notational convenience.

\begin{definition}
\label{h.b.def} Let $\mcC$ be a strict monoidal category and $x \in \mcC$.
A half-braiding $e_x$ for $x$ is a family of isomorphisms $\{e_x(y) \in \Hom_\mcC(xy,yx)\}_{y \in \mcC}$ satisfying
\begin{enumerate}

\item {Naturality:} for all $f \in \Hom(y,z), \quad (f \otimes \id_x) \circ e_x(y) = e_x(z) \circ (id_x\otimes f).$
\item {Half-braiding:} for all $y,z \in \mcC , \quad e_x(y\otimes z)= (id_y \otimes e_x(z)) \circ (e_x(y) \otimes \id_z).$
\item {Unit property:} $e_x(\unit)=\id_x$.

\end{enumerate}
\end{definition}
The objects in the Drinfeld center $\mcZ(\mcC)$ are direct sums of simple objects in $\mcC$ that admit half-braidings and the morphisms are those that behave compatibly with the half-braidings.  More precisely:
\begin{definition}
\label{qd} The Drinfeld center $\mcZ(\mcC)$ of a strict monoidal category
$\mcC$ has as objects pairs $(x,e_x)$, where $x \in \mcC$ and $e_x$ is a
half-braiding. The morphisms are given by \[\Hom \bigl((x,e_x),(y,e_y)\bigr)=
\bigl\{ f  \in \Hom_\mcC(x,y) \mid (\id_z \otimes f) \circ e_x(z) = e_y(z) \circ (f \otimes \id_z) \quad \forall z \in \mcC \bigr\}.\]
The tensor product of objects is given by $(x,e_x) \otimes (y,e_y)
= (xy,e_{xy})$, where
\[e_{xy}(z) = (e_x(z) \otimes \id_y) \circ (\id_x \otimes e_y(z)).\]
The tensor unit is $(\unit,e_\unit)$ where $e_\unit(x)=\id_x$. The composition
and tensor product of morphisms are inherited from $\mcC$. The
braiding is given by $c_{(x,e_x),(y,e_y)} = e_x(y)$.
\end{definition}
In \cite{MII} M\"uger proved the following:
\begin{theorem}
If $\mcC$ is a spherical fusion category, then $\mcZ(\mcC)$ is modular.
\end{theorem}

In the case that $\mcC$ is already modular the Drinfeld center factors $\mcZ(\mcC)=\mcC\boxtimes\mcC^{op}$ where $\mcC^{op}$ is the opposite category of $\mcC$ with the opposite braiding.  The topological central charge of $\mcZ(\mcC)$ is always $0$ mod $8$ because the resulting representations of mapping class groups are always linear.  Two spherical fusion categories $\mcD$ and $\mcB$ such that $\mcZ(\mcD)\cong \mcZ(\mcB)$ are called \emph{Morita equivalent}.

\subsubsection{Equivariantization, de-equivariantization, coring and gauging}

Finite groups can appear as certain symmetries of fusion categories, which can be exploited to produce new categories.  The easiest to understand is {de-equivariantization} \cite{MAlg,Brug}, which we will describe in the braided case.

Suppose that $\mcC$ is a ribbon fusion category with $\Rep(G)\cong\mcE\subset\mcC$ as a (symmetric, Tannakian) ribbon subcategory.  Then the algebra $\Gamma$ of functions on $G$ acts on $\mcC$ so that we may consider the category $\mcC_G$ of $\Gamma$-modules in $\mcC$, called the \emph{$G$-de-equivariantization} of $\mcC$.  The resulting category $\mcC_G$ is faithfully graded by $G$: $\mcC_G\cong \bigoplus_g (\mcC_G)_g$.

There are two particularly interesting cases to consider: 1) $\mcC$ is modular and 2)  $\mcE=\mcC^\prime$ the M\"uger center.  Notice that in case 1) $\mcC^\prime$ is trivial.
If $\mcC$ is modular, then the trivial component $(\mcC_G)_1$ is again modular \cite{DGNO1}, with $\dim((\mcC_G)_1)=\frac{\dim(\mcC)}{|G|^2}$.  In the physics literature passing from $\mcC$ to $(\mcC_G)_1$ is called \emph{boson condensation}, and the other components $(\mcC_G)_g$ are called the \emph{confined sectors} consisting of \emph{defects}.  In the special case that $\mcE$ is the maximal Tannakian symmetric subcategory then $(\mcC_G)_1$ is called the \emph{core} \cite{DGNO1}. For example $\Rep(G)$ is the maximal Tannakian subcategory of $\Rep(DG)$ and the core is just $\vect$.
If $\mcE=\mcC^\prime$ then $\mcC_G$ itself is modular, and it called the \emph{modularization} \cite{Brug}.  Modularization and coring may produce interesting new examples, but more typically these tools are used to reduce complicated categories to better understood, modular categories of smaller dimension: for example the $\mbbZ_2$-de-equivariantization of the categories $SO(2k+1)_2$ are TY whose trivial components are pointed modular categories.

The reverse process to boson condensation, called \emph{gauging} \cite{CGPW}, is somewhat more intricate, with both obstructions and choice dependency \cite{ENO3}.  One starts with a modular category $\mcC$ with an action of finite group $G$ via a homomorphism $\rho:G\rightarrow \Aut_\otimes^{br}(\mcC)$ to the (finite, see \cite{ENO3}) group of equivalence classes of braided tensor auto-equivalences of $\mcC$.  Producing a faithfully $G$-graded fusion category with trivial component $\mcC$--a so-called \emph{$G$-extension}--is the first hurdle: there are several potential (cohomological) obstructions, and even when the obstructions vanish there can be many such $G$-extensions.  The resulting fusion category $\mcC_{G}^{\times}$ is typically not modular (or even braided), instead it will be $G$-crossed modular category \cite{kirillov2004g}.  However $\mcC_{G}^{\times}$ comes with an action of $G$, and we may \emph{equivariantize} $\mcC_{G}^{\times}$--the inverse process to de-equivariantization--to obtain a modular category $(\mcC_{G}^{\times})^G$.  While verifying the obstructions vanish and parameterizing the choices can be difficult, we are often rewarded with interesting categories. To give an example, the pointed rank $4$, ``3-fermion theory'' $SO(8)_1$ has an action by $\mbbZ_3$ permuting the 3 fermions.  One $\mbbZ_3$-gauging produces the category $SU(3)_3$ which is an integral modular category of dimension $36$.

An important class of weakly integral categories are the \emph{weakly group theoretical} fusion categories, i.e. those that are Morita equivalent to a nilpotent category \cite{ENO2}.  Thanks to a recent paper of Natale \cite{Natcore} we have an alternative definition in the modular setting: a modular category $\mcC$ is \emph{weakly group theoretical} if the core of $\mcC$ is either pointed or the Deligne product of a pointed category with an Ising category (a non-integral modular category of dimension $4$).  In particular, any weakly group theoretical modular category is obtained by gauging Deligne products of Ising and pointed modular categories.

\begin{remark}
Can every modular category be constructed, using the above tools, from finite groups and quantum groups?
Most of the known modular categories do come from such constructions. However, exotic fusion categories arise in subfactor theory \cite{MN10}.  On one hand, taking the Drinfeld center of those unitary fusion categories gives rise to new unitary modular categories which do not resemble those constructed from quantum groups.  On the other hand, it is a folklore conjecture that all unitary modular categories can be generated from quantum groups \cite{HRW}.
  Loosely, we would like to call any unitary modular category that cannot be constructed from quantum group categories an exotic modular category.  But it is difficult to mathematically delineate all modular categories from quantum group constructions, thus define exoticness.  A first approximation using the Witt group \cite{DMNO} is described below, suggesting that, modulo Drinfeld centers, there are \emph{no} exotic modular categories.
\end{remark}

\subsection{Invariants of Modular Categories}

To classify and distinguish modular categories, we need invariants.  The theorem in {\cite{Tu94}} that each modular category $\mcC$ leads to a $(2+1)$-TQFT $(V,Z)$ can be regarded as a pairing between modular categories $\mcC$ and manifolds: $(\mcC, Y)=V(Y)$ regarded as a representation of the mapping class group for a $2$-manifold $Y$ or $(\mcC, X)=Z(X)$ the partition function for a $3$-manifold $X$ with some framing, possibly with a link $L$ inside.  The invariant $Z(X)$ is constructed using Kirby diagram of three manifolds.  Then for each fixed manifold $M$, we obtain invariants $(\mcC,M)$ of modular categories $\mcC$.  The most useful choices are the $2$-torus $T^2$ and some links in the $3$-sphere $S^3$.

\subsubsection{Modular data and $(S,T)$-uniqueness conjecture}

The most basic invariant of a modular category is the rank $|\Pi_\mcC|$, which can be realized as the dimension of the mapping class group representation of $T^2$.  Other useful invariants include the Grothendieck semiring $K_0(\mcC)$ (see section \ref{groth and dim}), the collection of invariants $\{d_a\}$ for the unknot
  colored by the label $a\in \Pi_\mcC$, i.e. the quantum dimensions and 
 the number $D=\sqrt{\sum_{a\in \Pi_\mcC}d_a^2}$.  The invariant of the Hopf
  link colored by $a,b$ are the entries of the $S$-matrix, see Figure \ref{Sab}. 
\begin{figure}[h]\includegraphics[width=1.8in]{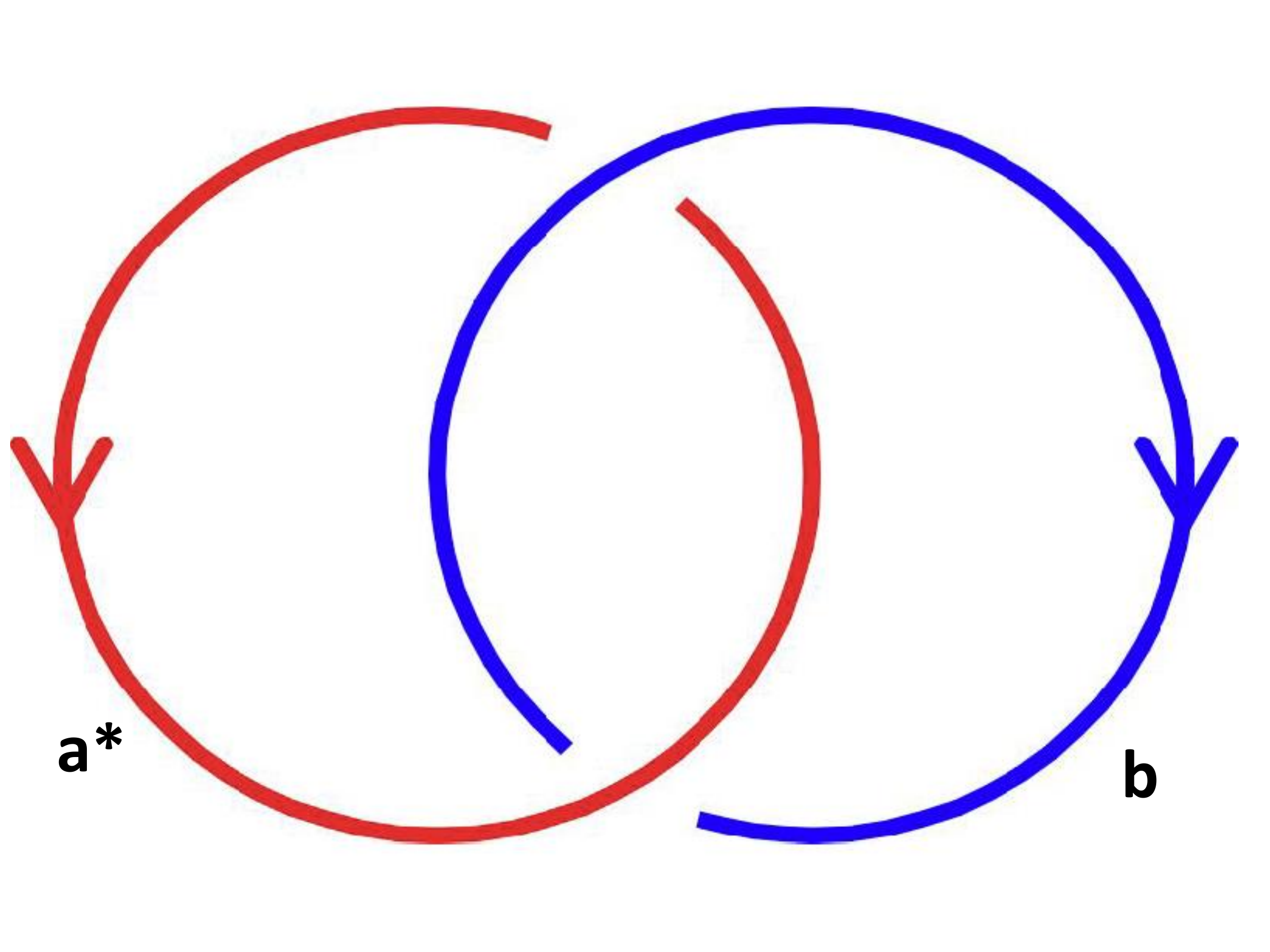}
\caption{\label{Sab}The $S$-matrix entry $S_{a,b}$}
\end{figure}
  
  The link invariant of the unknot with a right-handed
  kink colored by $a$ is $\theta_a\cdot d_a$ where $\theta_a$ is the {\it topological twist} of the label $a$.  The topological
  twists are encoded in a diagonal matrix $T=(\delta_{ab}\theta_a), a,b\in \Pi_\mcC$.   The
  $S$-matrix and $T$-matrix together define a projective representation of the
  modular group $\SL (2,\mathbb{Z})$ via (see \cite{Tu94, BKi}):
  \begin{equation*}
     \fs=\begin{pmatrix}0 & -1\\ 1 & 0 \end{pmatrix}\mapsto S,\qquad \ft=\begin{pmatrix} 1 & 1\\ 0 & 1 \end{pmatrix}\mapsto T.
  \end{equation*}

  Amazingly, the kernel of this projective representation of $\mcC$ is always a congruence
  subgroup of $\SL{(2,\mathbb{Z})}$, as was conjectured by Coste and Gannon \cite{CT94} and proved recently:

\begin{theorem}\cite{NS3}\label{congruence subgroup}
Let $\mcC$ be a modular category with $\ord(T)=N$ and
$\rho:(\mathfrak{s},\mathfrak{t})\rightarrow (S,T)$ the corresponding
$\SL{(2,\mathbb{Z})}$ representation.  Then $\ker\rho$ is a congruence subgroup of level $N$.\end{theorem}

The $S$-matrix determines the fusion rules through the Verlinde formula, and the $T$-matrix has finite order $\ord(T)$ by
  Vafa's theorem \cite{BKi}.  Together, the pair $(S,T)$ is called the \emph{modular data} of the category $\mcC$, and they satisfy many algebraic conditions \cite[Definition 2.7]{BNRW152}.  In particular the entries of $S$ and $T$ lie in the cyclotomic field $\mbbQ_N:=\mbbQ(e^{2\pi i/N})$ where  $\ord(T)=N$ and hence have abelian Galois groups.  Moreover, we have:
  \begin{theorem}\label{galgroupthm} Suppose $(S,T)$ is the modular data of a rank $r$ modular category $\mcC$ and let $\mbbQ(S)$ be the field generated over $\mbbQ$ by the entries of $S$.  \begin{enumerate}\item The Galois group $\Aut_{\mbbQ}(\mbbQ(S))$ is isomorphic to a subgroup of the symmetric group $\fS_r$.
  \item $\Aut_{\mbbQ(S)}(\mbbQ_N)\cong(\mbbZ/2\mbbZ)^\ell$.
  \end{enumerate}
  \end{theorem}
  
  The following is open:
  \begin{conjecture}
   The $S$ and $T$ matrices of a modular category $\mcC$ determine it up to ribbon equivalence.
  \end{conjecture}

Counterexamples to this conjecture have appeared recently \cite{ModularData17}.  

\subsubsection{Frobenius-Schur indicators}

For a finite group $G$, the $n$th Frobenius-Schur indicator of a representation $V$ over $\mbbC$ with character $\chi_V$ is given by $\nu_n(V):=\frac{1}{|G|}\sum_{g\in G}\chi_V(g^n)$. The classical Frobenius-Schur Theorem asserts that the second FS-indicator $\nu_2(V)$ of an irreducible representation $V$ must be 1, -1 or 0, which can be determined by the existence of non-degenerate $G$-invariant symmetric (or skew-symmetric) bilinear form on $V$. Moreover, the indicator value 1, -1 or 0 indicates respectively whether $V$ is real, pseudo-real or complex.

An expression for the second FS-indicator for each primary field of a rational conformal field theory were introduced by Bantay \cite{Ban1} in terms of modular data, which is more generally valid for simple objects in a modular category.  Higher FS-indicators for pivotal categories were developed by Ng and Schauenburg, in particular, for spherical fusion categories $\mcC$ over $\mbbC$ \cite{NS2}.  A general formula for the $n$th FS-indicator for simple objects in a modular category $\mcC$ is given by:
$$\nu_n(X_k): = \frac{1}{\dim{\mcC}} \sum_{i,j =0}^{r-1} N_{i,j}^k\,
          d_i d_j\left(\frac{\theta_i}{\theta_j}\right)^n.$$

For a spherical fusion category $\mcC$, the minimal positive integer $N$ such that $\nu_N(X_k)=d_k$ for all $k \in \Pi_\mcC$ is called the \emph{FS-exponent} of $\mcC$, and is denoted $\FSexp(\mcC)$ \cite{NS2}. For a modular category one can show that $\FSexp(\mcC)=\ord(T)$.  The FS-exponent of a spherical fusion category behaves, in many ways, like the exponent of a finite group. In fact, the FS-exponent of $Rep(G)$ for any finite group $G$ is equal to the exponent of $G$. These generalized FS-indicators for modular categories have surprising arithmetic properties, which play a key role in the proofs of several longstanding conjectures such as Theorem \ref{congruence subgroup}, the Galois Symmetry Theorem \cite[Theorem II.iii]{DLN1} and Theorem \ref{rank-finiteness}.

A crucial result relating the FS-exponent of a category to its dimension is the Cauchy theorem:
\begin{theorem}\cite{BNRW151}\label{Cauchy}
If $\mcC$ is a spherical fusion category over $\mbbC$, the set of prime ideals dividing the principal ideal $\langle\dim(\mcC)\rangle$ in the Dedekind domain $\mbbZ[e^{2\pi i/N}]$ is identical to that of $\langle N\rangle$ where $\FSexp(\mcC)=N$.
\end{theorem}
For the spherical category $\Rep(G)$ this is just a restatement of the classical Lagrange and Cauchy theorems, as the dimension is $|G|$ and the FS-exponent is the usual group exponent.
\begin{remark}
Each quantum group category $\mcC(\mathfrak{g},q,\ell)$ has a \emph{unimodal} version, i.e. with $\nu_2(V)=1$ for all self-dual simple objects, by choosing different ribbon elements \cite{Tu94}, but in general unitarity and modularity cannot both be preserved.  This can be easily seen from the semion theory $\mcC(\mathfrak{sl}_2,q,4)$.  The non-trivial simple object $s$ has quantum dimension=$1$, topological twist=$i$, and $\nu_2(s)=-1$.  If we choose a different ribbon element to obtain $\nu_2(s)=1$, then the resulting category is either non-unitary (the quantum dimension of $s$ would be =$-1$) or not modular (the topological twist of $s$ would be =$\pm 1$).   Remark 3 in XI.6.4 on page 512 of \cite{Tu94} leads to some confusion as the resulting unimodal categories are not in general the same as the quantum group types $\mcC(\mathfrak{g},q,\ell)$ for Reshetikhin-Turaev TQFTs.
\end{remark}

\subsection{Structure of Modular Categories}
A hypothetical periodic table for modular categories should be grouped by families of related categories, e.g. \cite{LW17}. Unfortunately, a satisfactory structure theorem for modular categories is still lacking.

\subsubsection{Prime decomposition}
If $\mcB\subset\mcC$ are both modular categories, then $\mcC\cong\mcB\boxtimes \mcD$ for some modular category $\mcD$ \cite{M2}.  If $\mcC\not\cong\vect$ contains no modular subcategories then $\mcC$ is called \emph{prime} and M\"uger proved the following prime decomposition theorem:
\begin{theorem}\cite{M2}
 Every modular category is a product of prime modular categories.
\end{theorem}

For example, for $\gcd(N,k)=1$ we have a factorization as modular categories $SU(N)_k\cong SU(N)_1\boxtimes PSU(N)_k$, where $SU(N)_1$ is a pointed modular category with fusion rules like the group $\mbbZ/N\mbbZ$.

Unfortunately this decomposition is not unique: for example, the squares of pointed modular categories can coincide: $SO(16)_1^{\boxtimes 2}\cong SO(8)_1^{\boxtimes 2}$.  However, if there are no non-trivial $1$-dimensional objects the prime decomposition is unique \cite{M2}.  

\subsubsection{The Modular Witt group}

The classical Witt group of quadratic forms on finite abelian groups is important for many applications such as surgery theory in topology.
Recently a similar theory has been under development for nondegenerate braided fusion categories \cite{DMNO,DNO} and generalized further to braided fusion categories over symmetric fusion categories \cite{DNO}, and we give a flavor of the theory for modular (i.e. spherical nondegenerate braided fusion) categories.
\begin{definition}
Two modular categories $\mcC_1,\mcC_2$ are Witt equivalent if there exist spherical fusion categories $A_1$ and $A_2$ such that
$ \mcC_1 \boxtimes \mcZ(A_1) \simeq \mcC_2 \boxtimes \mcZ(A_2)$
where $\mcZ(A_i)$ are Drinfeld centers, and $\simeq$ is ribbon equivalence.
\end{definition}
Witt equivalence is an equivalence relation, and $\boxtimes$ descends to Witt classes.
\begin{theorem}\cite{DMNO}
Witt classes form an abelian group under $\boxtimes$, and $\mcC$ is in the trivial class if and only if $\mcC \simeq \mcZ(A)$ for some spherical fusion category $A$.
\end{theorem}

The full structure of the Witt group $\mathcal{W}_{un}$ of UMCs is unknown.  It is known to be an infinite group and the torsion subgroup is a 2-group, and the maximal finite order of an element is $32$ \cite{DNO}.  One application of the Witt group is that it provides a precise (if somewhat coarse) framework for studying exoticness: Question 6.4 of \cite{DMNO} asks if the Witt group $\mcW_{un}$ of unitary modular categories is generated by the classes of quantum groups.
There is an obvious homomorphism $\tilde{c}_{\mathrm{top}} \colon \mathcal{W}_{un} \to \frac{\mbbQ}{8\mbbZ}$.
Are there  nontrivial homomorphisms other than the one given by the  topological central charge?

\subsubsection{Integrality and group-like properties}
Known examples of weakly integral categories have several distinguishing characteristics, which suggests that a  structural description of this subclass is within reach.  One example is in the Property $F$ conjecture \cite{ERW,NR11}:

\begin{conjecture}\label{propF} Let $X,Y\in\mcC$ be a simple objects in a braided fusion category. The associated braid group representation $\rho_n^X(Y):\mcB_n\rightarrow \Aut(\Hom(Y,X^{\otimes n}))$ has finite image if, and only if $\FPdim(X)^2\in \mbbZ$.
\end{conjecture}
In particular, the braid group representations associated with objects in weakly integral modular categories would have finite image which would imply that weakly integral modular categories model non-universal anyons.
There is significant empirical evidence: it is known to be true for group-theoretical categories (e.g. $\Rep(D^\omega G)$ see \cite{ERW}) and for quantum group categories (see e.g. \cite{FLW02,RWenzl,Ruma}).  More generally one can ask if \emph{all} mapping class group representations obtained from a weakly integral modular category have finite image, and this was recently answered in the affirmative for $\Rep(D^\omega G)$ \cite{Fuchs15,Gustafson}.

Every weakly group-theoretical modular category is weakly integral, and  the converse is conjectured in \cite{ENO2}.  In particular this would imply that every weakly integral modular category is obtained from products of Ising and pointed modular categories by gauging.  If true, this suggests a route to a proof of Conjecture \ref{propF}: relate the braid group representations coming from a weakly integral $\mcC$ to those of its core $(\mcC_G)_1$.

In \cite{twoparas} it is suggested that weakly integral modular categories correspond to link invariants that are approximable in polynomial time. Yet another potential characterization of weakly integral categories is the following, which has been verified for the Jones representations of $\mcB_n$:
\begin{conjecture}\label{localizeconj}\cite{RWlocal} A simple object $X\in\mcC$ in a unitary braided fusion category $\mcC$ can be localized (see Definition \ref{localize}) if, and only if, $\FPdim(X)^2\in\mbbZ$.
\end{conjecture}

It is not known if every integral modular category is weakly group-theoretical.  An open problem posed by Nikshych\footnote{Simons Center, September 2015} is to show that any non-pointed integral modular category contains a nontrivial Tannakian (or even just symmetric) subcategory.  

\subsubsection{Rank-finiteness and low-rank classification}
Combining Theorems \ref{galgroupthm} and \ref{Cauchy} with some results in analytic number theory \cite{evertse}, we obtain the rank-finiteness theorem, originally conjectured by the second author in 2003:
\begin{theorem}\label{rank-finiteness}\cite{BNRW151}
There are only finitely many modular categories of rank $r$, up to equivalence.
\end{theorem}
This result shows that a program classifying modular categories by rank is, in principle, possible.  Although Etingof has shown \cite[Remark 4.5]{BNRW151} that the number of weakly integral rank $r$ modular categories grows faster than any polynomial in $r$, excising them might leave a more tractable class to consider.

Modular categories of rank$\leq 5$ have been completely classified \cite{RSW09,BNRW152} up to braided monoidal equivalence.  The following list contains one representative from each Grothendieck equivalence class of prime modular categories of low rank (see \cite{RSW09,BNRW152}):
\begin{enumerate}
\item[2.] Rank=2: $PSU(2)_3$ (Fibonacci), $SU(2)_1$ (Semion, pointed)
\item[3.] Rank=3: $SU(3)_1$ (pointed), $SU(2)_2$ (Ising, weakly integral), $PSU(2)_5$
\item[4.] Rank=4: $SU(4)_1$ (pointed), $SO(8)_1$ (pointed), $PSU(2)_7$
\item[5.] Rank=5: $SU(5)_1$ (pointed), $SU(2)_4$ (weakly integral), $PSU(2)_9$, $PSU(3)_4$.
\end{enumerate}

\section{Extensions and More Open Problems}

TQC with anyons is relatively mature, but has many extensions to TQC with symmetry defects, gapped boundaries and the defects between them, and extended objects in higher dimensions.  Firstly, our discussion of TQC so far is based on TPMs of boson systems, but real topological materials such as the fractional quantum Hall liquids are fermion systems.  Therefore, we need a theory for fermionic TPMs.  Secondly, topology and conventional group symmetry have interesting interplay as illustrated by topological insulators and superconductors. TPMs with group symmetry can support symmetry defects.   Thirdly, real samples have boundaries, so the boundary physics and the correspondence with the bulk (interior) is also very rich.  Lastly, while it is possible to engineer two dimensional TPMs, three dimensional materials are much more common.  We discuss these extensions in this section.  More speculative extensions can be found \cite{WangBA}.

\subsection{Fermions}

The most important class of TPMs is two dimensional electron liquids which exhibit the fractional quantum Hall effect (see \cite{NSSFD08} and references therein).  Usually fractional quantum Hall liquids are modeled by Witten-Chern-Simons TQFTs at low energy based on bosonization such as flux attachment.  But subtle effects due to the fermionic nature of electrons are better modeled by refined theories of TQFTs (or UMCs) such as spin TQFTs (or fermionic modular categories) \cite{BM05}.
A refinement of unitary modular categories to spin modular categories and their local  sectors---super-modular categories has been studied \cite{BG17}.

\subsubsection{Spin TQFTs}

\begin{definition}

A \emph{spin modular category} is a unitary modular category $\mcB$ with a chosen invertible object $f$ with $\theta_f=-1$.  An invertible object $f$ with $\theta_f=-1$ is called a fermion.

\end{definition}

Let $Bord_{2,3}^{\textrm{spin}}$ be the spin bordism category of spin $2$- and $3$-manifolds.  The objects $(Y,\sigma)$ of $Bord_{2,3}^{\textrm{spin}}$ are oriented surfaces $Y$ with spin structures $\sigma$ (a lifting of the $SO$-frame bundle to a $Spin$-frame bundle), and morphisms are equivalence classes of spin-bordisms.  Let s-$Vec$ be the category of super vector spaces and even linear maps, which is a symmetric fusion category.

\begin{definition}

A spin TQFT is a symmetric monoidal projective functor $(V^s,Z^s)$ from $Bord_{2,3}^{\textrm{spin}}$ to the symmetric fusion category s-$Vec$.

\end{definition}

Each spin modular category $(\mcB, f)$ gives rise to a spin TQFT by decomposing the TQFT associated to the spin modular category regarded just as a modular category as follows.

Given a spin modular category $(\mcB, f)$, then there is a TQFT $(V,Z)$ as constructed in \cite{Tu94} from the modular category $\mcB$.  The partition function $Z(X^3)$ of an oriented closed $3$-manifold $X$ will be decomposed as a sum $Z(X^3)=\sum_\sigma Z(X^3,\sigma)$, where $\sigma$ is a spin structure of $X$.  Hence $Z(X^3,\sigma)$ is an invariant for spin closed oriented $3$-manifolds. For each oriented closed surface $Y$, the TQFT Hilbert space $V_{\mcB}(Y)$ is decomposed into subspaces indexed by the spin structures of $Y$: $V_\mcB(Y)=\oplus_\sigma V(Y,\sigma)$.

For simplicity, we will only define the Hilbert space $V^s(Y,\sigma)$ for a spin closed surface $(Y,\sigma)$.  Set $V^s_0(Y,\sigma)=V(Y,\sigma)$ and $V^s_1(Y,\sigma)=V(Y_0,f;\sigma)$, where $V(Y_0,f;\sigma)$ is the Hilbert space associated to the punctured spin surface $Y_0$ with a single puncture of $Y$ labeled by the fermion $f$.  Then setting $Z^s(X^3,\sigma)$ equal to the invariant from the decomposition of $Z(X^3)$, and $V^s(Y,\sigma)=V^s_0\oplus V^s_1$ leads to a spin TQFT $(V^s,Z^s)$.

It is easy to check that while the disjoint union axiom does not hold, the $\mathbb{Z}_2$ version of the disjoint union axiom does hold.

\subsubsection{16-fold way}

\begin{definition}

A super-modular category is a unitary pre-modular category $\mcB$ whose M\"uger center is isomorphic to $\sVec$, the symmetric fusion category $\sVec$ generated by two simple objects $\{1,f\}$ for some fermion $f$.

\end{definition}

Fermion systems have a fermion number operator $(-1)^F$ which leads to the fermion parity: eigenstates of $(-1)^F$ with eigenvalue $+1$ are states with an even number of fermions and  eigenstates of $(-1)^F$ with eigenvalue $-1$ are states with an odd number of fermions.  This fermion parity is like a $\mathbb{Z}_2$-symmetry in many ways, but it is not strictly a symmetry because fermion parity cannot be broken.  Nevertheless, we can consider the gauging of the fermion parity (compare with \cite{BBCW,CGPW}).  In our model, the gaugings of the fermion parity are the minimal extensions of the super-modular category $\mcB$ to its covering spin modular categories $\mcC$.  In two spatial dimensions, gauging the fermion parity seems to be un-obstructed.  So we conjecture that a minimal modular extension always exists, and there are exactly $16$ such minimal extensions of super-modular categories.  We will refer to this as the $16$-fold way conjecture \cite{ENO3,BG17}:

\begin{conjecture}\label{16foldway}
Let $\mcB$ be super-modular.  Then $\mcB$ has precisely $16$ minimal unitary modular extensions.
\end{conjecture}
A topological approach to this conjecture would be to construct a spin TQFT for each super-modular category.  The $16$ here is probably the same as in Rochlin's theorem.
It is known that if $\mcB$ has one minimal modular extensions then it has precisely 16 \cite{LKW17}.

\subsection{Symmetry Defects}

\subsubsection{Representation of Modular Categories}
As modular categories are categorifications of rings, module categories categorify representations of rings.
 Suppose $\mcC$ is a modular category, a \emph{module category} $\mathcal{M}$ over $\mcC$ is a categorical representation of $\mcC$.  A \emph{left module category} $\mcM$ over $\mcC$ is a semi-simple category with a bi-functor $\alpha_{\mcM}: \mcC\times \mcM \rightarrow \mcM$ that satisfies the analogues of pentagons and the unit axiom.   Right module categories are similarly defined, and a bi-module category is a simultaneously left and right module category such that the left and right actions are compatible.  Bi-module categories can be tensored together just like bi-modules over algebras.  Since $\mcC$ is braided, a left module category is naturally a bi-module category using the braiding.  A bi-module category $\mcM$ over $\mcC$ is invertible if there is another bi-module category $\mathcal{N}$ such that $\mcM \boxtimes \mathcal{N}$ and $\mathcal{N} \boxtimes \mcM$ are both equivalent to $\mcC$---the trivial bi-module category over $\mcC$.  The invertible module categories over $\mcC$ form the Picard categorical-group $\underline{\textrm{Pic}}(\mcC)$.  As $\mcC$ is a modular category there is an isomorphism of categorical-groups ${\uPic}(\mcC)\cong \underline{\Aut_\otimes^{br}}(\mcC)$ \cite{ENO2}.
 
Th first obstruction to gauging an action $\rho:G\rightarrow \Aut_\otimes^{br}(\mcC)$ on a UMC $\mcC$ is to lift it to a topological symmetry:
\begin{definition}
A finite group $G$ is a \emph{topological symmetry} of a UMC $\mcC$ if there is a monoidal functor $\underline{\rho}: \underline{G}\rightarrow  \underline{\Aut_\otimes^{br}}(\mcC)\cong\uPic(\mcC)$, where $\underline{G}$ is the categorical group with a single object and every group element is an invertible morphism.
%We will denote the topological symmetry as $(\underline{\rho}, G)$ or simply %$\underline{\rho}$ and say that $G$ acts categorically on $\mcC$.
\end{definition}
The resulting $G$-extension is a $G$-crossed braided fusion category, therefore the action is by conjugation and it follows that if a defect is fixed, then the group element is in the center.

There are two obstructions---one to the existence of tensor product and the other to the associativity.  
If all obstructions vanish, then the $G$-extension step of gauging is to add (extrinsic topological) defects $X_g$, which are objects in the invertible module category $\mcC_g$ corresponding to $g$ under the above functor.

If a defect $X_g$ is fixed by the $G$ action ($g$ has to be in the center of $G$), then under equivariantization, the defect $X_g$ becomes several anyons $(X_g, \pi)$, where $\pi$ is some irreducible representation of $G$.  Let $\rho_{X_g, n}$ be the projective representation of the braid group $\mcB_n$ from the $G$-crossing of $X_g$, and $\rho_{(X_g, \pi),n}$ be the representation of $\mcB_n$ from the anyon $(X_g, \pi)$ in the gauged modular category, where $\pi$ is a projective representation.

\begin{conjecture}

$\rho_{X_g, n}$ is equivalent to $\rho_{(X_g, \pi),n}$ as projective representations for any $n$ and $\pi$, where $X_g$ is fixed by the $G$ action.

\end{conjecture}

Such projective representations of the braid group from symmetry defects can be used for quantum computing and enhance the computational power of anyons \cite{delaney17SD}.  In the case of bilayer Ising theory, the Ising anyon $\sigma$ can be made universal using symmetric defect states as ancillas \cite{BF16}.

\subsection{Boundaries}

How to model the boundary physics of a TPM is still a subtle question.

\subsubsection{Gapped Boundaries}

Recent studies of TPMs revealed that certain TPMs also support gapped boundaries \cite{Bravyi98}.
In the UMC model of a 2D doubled topological order $\mathcal{B}=\mathcal{Z}(\mathcal{C})$, a stable gapped boundary or gapped hole is modeled by a Lagrangian algebra $\mathcal{A}$ in $\mathcal{B}$.  The Lagrangian algebra $\mathcal{A}$ consists of a collection of bulk bosonic anyons that can be condensed to vacuum at the boundary, and the corresponding gapped boundary is a condensate of those anyons which behaves as a non-abelian anyon of quantum dimension $d_{\mathcal{A}}$.  Lagrangian algebras in $\mathcal{B}=\mathcal{Z}(\mathcal{C})$ are in one-to-one correspondence with indecomposable module categories $\mathcal{M}$  over $\mathcal{C}$, which can also be used to label gapped boundaries and used for TQC \cite{KK12,CongCMP17,CMW16,CongPRL17,cong2017survey}.

If a single boundary circle is divided into many segments and each segment labeled by an indecomposable module $\mcM_i$, then the defects between different boundary segments are modeled by functors $Fun_{\mcB}(\mcM_i,\mcM_j)$ \cite{KK12, CMW16,CMW17}.  These boundary defects can support degeneracy and projective representations of the braid group.  In certain cases, they are related to symmetry defects and it is conjectured their corresponding braid group representations are projectively equivalent.
They can also be used for TQC \cite{CMW17}.

\subsubsection{Gapless Boundaries}

It is widely believed in the case of fractional quantum Hall states that the boundary physics can be modeled by a unitary chiral CFTs ($\chi$CFT) $\mcV$ \cite{Wen92,Read09}.
As an instance of a bulk-edge correspondence, the UMC $\mcC_\mcV$ encoded in the boundary $\chi$CFT $\mcV$ is the same as the UMC $\mcC_\mcH$ of the bulk.
Moreover, the UMC $\mcC_\mcH$ has a multiplicative central charge $\chi=e^{\pi ic/4}$, where $c$, called the (chiral) topological central charge of $\mcC_\mcH$, is a non-negative rational number defined modulo $8$, which
agrees with the central charge of $\mcV$.

It is conjectured that this bulk-edge correspondence exists beyond the fractional quantum Hall states, so that for any given UMC $\mcC$, there is always a unitary $\chi$CFT $\mathcal{V}$ such that its UMC $\mcC_\mcV$ is $\mcC$ and its central charge is equal to the topological central charge of $\mcC$, modulo $8$ \cite{TW17}.  The same conjecture was made by Gannon \footnote{Casa Matematica Oaxaca, August 2016} as an analogue of Tannaka-Krein duality.
A long-term goal is to classify unitary $\chi$CFTs based on progress in classifying UMCs \cite{BNRW152,RSW09}.

Bulk-edge correspondence is the topological analogue of Ads/CFT with topological phase replacing quantum gravity.  But the detailed correspondence can be very subtle.  In the mathematical context, it is the question how to construct a VOA from its representation modular category.

\subsubsection{From UMCs to $\chi$CFTs}

We only consider unitary $\chi$CFTs and will use VOAs as our mathematical $\chi$CFT.
The minimal energies (or minimal conformal weights) $\{h_i\}$ of a nice VOA are encoded, mod 1, in the exponents of the topological twists of its UMC by $\theta_i = e^{2 \pi i h_i}$.
Therefore, one set of natural extra data to consider would be a lifting of the exponents of the topological twists from equivalence classes of rational numbers (modulo $1$) to actual rational numbers.
Since we are mainly interested in unitary theories, we only consider liftings for which $h_i \ge 0$ for all $i$.
It is not impossible that a consistent lifting of the topological twists is sufficient to determine a corresponding CFT within a given genus, at least when the CFT has non-trivial
representation theory (i.e. when the CFT is not holomorphic).

\begin{conjecture}\label{cnjGenusFiniteness}
Given a UMC $\mcC$, 
there is a central charge $c$ such that the admissible genus $(\mcC, c)$ is realizable.
\end{conjecture}

\subsection{(3+1)-TQFTs and  $3D$ Topological Phases of Matter}

The most interesting future direction is in $(3+1)$-TQFTs.  Mathematically, $(3+1)$-TQFT that can distinguish smooth structures would be highly desirable for the classification of smooth $4$-manifolds.  Physically, $3D$ space is the real physical dimension.  One lesson we learn from lower dimensions is that we might also want to consider $(4+1)$-TQFTs because they would provide understanding of anomalous $(3+1)$-TQFTs.

\subsubsection{(3+1)-TQFTs from $G$-crossed categories and spherical $2$-fusion categories}

The most general construction so far for state sum $(3+1)$-TQFTs is based on unitary $G$-crossed braided fusion categories, which are special cases of the unknown spherical $2$-fusion categories \cite{Cui16}.  Lattice realization of these state sum $(3+1)$-TQFTs as TPMs is given in \cite{WW12,WW17}.  Spherical $2$-fusion categories should be the fully dualizable objects in the $4$-category target of fully extended $(3+1)$-TQFTs \cite{Lurie09}.

\begin{conjecture}

The partition function of any unitary $(3+1)$-TQFT is a homotopy invariant.

\end{conjecture}

A stronger version of the conjecture would be such TQFTs are simply generalizations of the Dijkgraaf-Witten TQFTs so that the partition functions count homotopy classes of maps between higher homotopy types.

\subsubsection{Representation of motion groups}

Any $(3+1)$-TQFT will provide representations of motion groups of any link $L$ in any $3$-manifold $Y$.  Very little is known about these representations, but the ubiquity of braid groups in 2D TQC models hints at a similar role for these representations.

The simplest motion group is that of the unlink of $n$ circles in $S^3$.
This motion group is generated by ``leapfrogging" the $i$th circle through the $(i+1)$st $\sigma_i$ and loop interchanges $s_i$ for $1\leq i\leq n-1$, and called the \textit{Loop Braid Group}, $\mathcal{LB}_n$. Abstractly, $\mcL\mcB_n$ is obtained from the free product $\mcB_n*\fS_n$ of the $n$-strand braid group generated by the $\sigma_i$ and the symmetric group generated by the $s_i$ by adding the
(mixed) relations:
\begin{equation*} \sigma_i\sigma_{i+1}s_{i}=
s_{i+1}\sigma_{i}\sigma_{i+1},\quad
 s_is_{i+1}\sigma_{i}=
\sigma_{i+1}s_{i}s_{i+1},\quad 1\leq i\leq n-2, \quad
 \sigma_is_j=s_j\sigma_i \quad
\text{if} \quad |i-j|>1.
\end{equation*}

This is a relatively new area of development, for which many questions and research directions remain unexplored.  A first mathematical step is to study the unitary representations of the loop braid group, which is already underway \cite{KMR17,MRC15}.  It is reasonable to consider other configurations, such as loops bound concentrically to an auxiliary loop, or knotted loops.

\section{Quantum Matters}

The most rigorous creation of the human mind is the mathematical world.  Equally impressive is our creation of the computing world. At this writing, Machine is beating the best GO player in the world.  Man, Machine, and Nature meet at TQC.  What will be the implication of TQC, if any, for the future of mathematics?  Two interpretations of quantum mathematics would be mathematics inspired by quantum principles or mathematics based on an unknown quantum logic.

\subsection{Quantum Logic}

Logic seems to be empirical, then would quantum mechanics change logic?  There are interesting research in quantum logics, and quantum information provides another reason to return to this issue \cite{DMW13,Aaronson13}.

\subsection{Complexity Classes as Mathematical Axioms}

Another direction from TQC is Freedman's suggestion of complexity classes as mathematical axioms.  Some interesting implications in topology from complexity theories can be found in \cite{Freedman09,CFW16}.

%\section*{Appendix of Notations and Terminologies}

%\bibliographystyle{amsplain}

%\bibliography{TqcReferences}

%\printbibliography
\bibliographystyle{abbrv}
\bibliography{TqcReferences.bib}
\end{document}